\documentclass[12pt,a4paper,fleqn]{article}
\usepackage{lscape,exscale,amsthm}
\usepackage[intlimits]{amsmath}
\usepackage{rawfonts}
\usepackage{latexsym}
\usepackage[cp850]{inputenc}
\usepackage{epsfig}
\usepackage{bbm}
\usepackage{enumitem}
\usepackage{float}

\usepackage{natbib}
\usepackage{blindtext}

\RequirePackage[OT1]{fontenc}
\RequirePackage{amsthm,amsmath}
\newcommand{\bSigma}{\mathbf{\Sigma}}

\newcommand{\bx}{\mathbf{x}}
\newcommand{\bX}{\mathbf{X}}

\newcommand{\by}{\mathbf{y}}

\newcommand{\bC}{\mathbf{C}}

\newcommand{\br}{\mathbf{r}}
\newcommand{\ba}{\mathbf{a}}

\newcommand{\bR}{\mathbf{R}}

\newcommand{\bz}{\mathbf{z}}
\newcommand{\bB}{\mathbf{B}}

\newcommand{\bi}{\mathbf{1}}

\newcommand{\bI}{\mathbf{I}}

\newcommand{\tr}{\operatorname{tr}}

\newcommand{\bA}{\mathbf{A}}

\newcommand{\bOmega}{\mathbf{\Omega}}

\newcommand{\bS}{\mathbf{S}}

\newcommand{\bZ}{\mathbf{Z}}

\usepackage{tikz}

\usepackage{stackrel}

\numberwithin{equation}{section}
\theoremstyle{plain}
\newtheorem{theorem}{Theorem}[section]

\newtheorem{lemma}{Lemma}[section]

\newtheorem{remark}{Remark}
\usepackage{color}

\usepackage{ulem}

\usepackage{hyperref}
\hypersetup{
urlcolor=black,
  menucolor=black,
  citecolor=black,
  anchorcolor=black,
  filecolor=black,
  linkcolor=black,
  colorlinks=true,
}
\usepackage{dsfont}
\newcommand{\E}{\mathbbm{E}}

\definecolor{darkblue}{rgb}{.1, 0.1,.8}
\definecolor{darkgreen}{rgb}{0,0.8,0.2}
\definecolor{darkred}{rgb}{.8, .1,.1}

\newcommand{\diag}{\operatorname{diag}}
\renewcommand{\P}{\mathbbm{P}}

\newcommand{\cid}{\stackrel{d}{\rightarrow}}
\newcommand{\cip}{\stackrel{\P}{\rightarrow}}
\newcommand{\norm}[1]{\|#1\|}
\newcommand{\nto}{n \to \infty}
\usepackage{slashbox}
\newcommand{\Var}{\operatorname{Var}}
\newcommand{\Cov}{\operatorname{Cov}}
\newcommand{\Corr}{\operatorname{Corr}}
\newcommand{\ov}{\overline}

\usepackage{bbm}
\begin{document}

%\today

\begin{center}
  \vspace*{2cm} \noindent {\bf \large Logarithmic law of large random correlation matrix}\\

  \vspace{1cm} \noindent {\sc\small Nestor Parolya$^a$\footnote{Corresponding author n.parolya@tudelft.nl}, Johannes Heiny$^b$ and Dorota Kurowicka$^a$}\\
\vspace{1cm}
{\it \footnotesize  $^a$
Delft Institute of Applied Mathematics, Delft University of Technology,
Delft, The Netherlands}\\
{\it \footnotesize  $^b$ Ruhr-University Bochum, Department of Mathematics, Bochum, Germany
}
 \end{center}
\vspace{1cm}

\begin{abstract}
Consider a random vector $\mathbf{y}=\mathbf{\Sigma}^{1/2}\mathbf{x}$, where the $p$ elements of the vector $\mathbf{x}$ are i.i.d. real-valued random variables with zero mean and finite fourth moment, and $\mathbf{\Sigma}^{1/2}$ is a deterministic $p\times p$ matrix such that the spectral norm of the population correlation matrix $\mathbf{R}$ of $\mathbf{y}$ is uniformly bounded. In this paper, we find that the log determinant of the sample correlation matrix $\hat{\mathbf{R}}$ based on a sample of size $n$ from the distribution of $\mathbf{y}$ satisfies a CLT (central limit theorem) for $p/n\to \gamma\in (0, 1]$ and $p\leq n$. Explicit formulas for the asymptotic mean and variance are provided.
 In case the mean of $\mathbf{y}$ is unknown, we show that after recentering by the empirical mean the obtained CLT holds with a shift in the asymptotic mean. This result is of independent interest in both large dimensional random matrix theory and high-dimensional statistical literature of large sample correlation matrices for non-normal data. At last, the obtained findings are applied for testing of uncorrelatedness of $p$ random variables. Surprisingly, in the null case $\mathbf{R}=\mathbf{I}$, the test statistic becomes completely pivotal and the extensive simulations show that the obtained CLT also holds if the moments of order four do not exist at all, which conjectures a promising and robust test statistic for heavy-tailed high-dimensional data.
\end{abstract}

\vspace{0.7cm}

 \noindent AMS 2010 subject classifications: 60B20, 60F05, 60F15, 60F17, 62H10\\
 \noindent {\it Keywords}: Sample correlation matrix, CLT, log determinant, large-dimensional asymptotics, random matrix theory, dependent data. \\
\section{Introduction}

Sample correlation matrices have always been of vital importance from both theoretical and practical points of view. Principal component analyis, for example, extracts valuable information about large data sets from the eigenvalues of the sample correlation matrix.
 In particular, the determinant of the sample correlation matrix is one of the most fundamental matrix functions and has been extensively studied in the theory of random matrices as well as in multivariate statistics (see, for instance, the classical monographs \cite{muirhead1982} and \cite{anderson2003}).
The determinant of a random correlation matrix has numerous applications in stochastic geometry (volume of parallelotope, see \cite{nielsen1999}) and hypothesis testing (likelihood ratio test) in multivariate statistics (see, \cite{anderson2003}). Its properties have been studied by many authors under various settings (see, e.g., \cite{nguyen2014} and references therein). For instance,
\cite{goodman1963} proved the central limit theorem (CLT) of the logarithmic determinant for random Gaussian matrices, \cite{tao2012} for Wigner matrices, \cite{nguyen2014} for  real i.i.d. random matrices under subexponential tail conditions, and \citet{bao2015, wang2018} for general i.i.d. matrices under existence of the 4th moments of matrix entries, to mention a few. In some special cases where a suitable stochastic representation is available, \cite{grote:kabluchko:thaele:2019} also proved large deviation results and \cite{heiny:johnston:prochno:2021} fast Berry--Esseen bounds. 

Our particular interest covers the sample correlation matrix denoted by $\hat\bR$ which is computed from a random sample $\by_1\ldots,\by_n$.
The determinant $\det\hat{\bR}$ is the likelihood ratio test statistic for testing the independence of the entries of a large dimensional random vector coming from a multivariate normal population. If the population correlation matrix is equal to identity, i.e., $\bR=\bI$, several results are available about its sample counterpart $\hat{\bR}$. In particular, under multivariate normality the density of $\det\hat{\bR}$ is proportional to $(\det R)^{(n-p-2)/2}dR$, see, \cite[Theorem 5.1.3]{muirhead1982}. However, the density of the eigenvalues of $\hat{\bR}$ cannot be obtained in the closed form, which makes the analysis of this random matrix challenging. Nevertheless, some asymptotic properties of the large sample correlation matrix $\hat\bR$  have been obtained in case of $\bR=\bI$ and normally distributed data. For example, the empirical distribution of eigenvalues of $\hat\bR$ follows the well-known Marchenko-Pastur law as first shown in \cite{jiang2004b} (see also \cite{bai:zhou:2008,heiny2018} for more general conditions) , while its largest eigenvalue obeys the Tracy-Widom law \cite{bao2012}. The properly normalized largest off-diagonal entry of $\hat{\bR}$ congerges to a Gumbel distribution as shown in \cite{jiang2004a} and later generalized in various directions in \cite{zhou:2007,li:qi:rosalsky:2012,liu:lin:shao:2008,cai:jiang:2011} and most recenty to a point process setting in \cite{heiny:mikosch:yslas:2020}.

 Moreover, under multivariate normality with $\bR=\bI$, the CLT for $\log\det\hat\bR$ holds (see, \citet{jiang2013, jiang2015}). These results were further generalized to non-normal populations by \cite{gao2017}. On the other hand, not much is known in case $\bR\neq\bI$. Recently the paper by \cite{jiang2019} sheds some light on this challenging case when the data comes from multivariate normal distribution. Particularly, \cite{jiang2019} showed that the properly normalized logarithmic determinant of $\hat\bR$ satisfies the CLT under some conditions on the eigenvalues of the population correlation matrix (minimum eigenvalue greater than $1/2$).

In this paper we will prove the CLT for $\log\det\hat{\bR}$ under very generic conditions on the data generating process, i.e., non-normality and general population correlation matrix $\bR\neq\bI$. We provide a closed form expression of the asymptotic mean and variance and discuss how these correspond to \cite{jiang2019} in case of normality. The work is under setting of increasing dimension $p$ and sample size $n$ diverging to infinity simultaneously, while their ratio tends to a constant $\gamma\in(0,1]$. Our results are  further applied to testing the uncorrelatedness of the elements of a high-dimensional random vector from arbitrary population with finite 4th moments.
Moreover, we investigate in detail the behavior of the log-determinant near singularity, i.e., $p/n\to1$ as $n\to\infty$, constructing of test on uniformity of entries of the sample correlation matrix. Interestingly, the distribution of the test statistic is independent of the population's fourth moments under $\bR=\bI$ in case of non-normal data. This property indicates a very well-behaved statistic for heavy-tailed distributions of the entries of large data matrix and opens a new direction for research in this topic.

Our paper is structured as follows: in Section~\ref{sec:2} we formulate the main result, namely the CLT for the logarithmic determinant of the sample correlation matrix for observations with mean zero. We also prove a similar CLT for observations with non-zero mean, for which an additional centering by the sample mean is needed. Section~\ref{sec:3} is devoted to applications of the obtained results. Here we study the behavior of the test statistic for testing the uncorrelatedness of the entries of large random vector and testing the uniformity of the entries of a large random correlation matrix. In Section~\ref{proofs} the proofs are given and Section~\ref{sec:5} with auxiliary lemmas finishes the paper.

\section{Logarithmic law of sample correlation matrix}\label{sec:2}

We consider a $p$-dimensional population $\by=\bSigma^{1/2}\bx$, where the $p$ elements of the vector $\bx$ are i.i.d. real-valued random variables and $\bSigma^{1/2}$ is a deterministic $p\times p$ matrix. The corresponding population correlation matrix of $\by$ is then given by $\bR = \diag(\bSigma)^{-1/2} \,\bSigma \, \diag(\bSigma)^{-1/2}$, where $\diag(\bSigma)$ denotes the diagonal matrix with the same diagonal elements as $\bSigma$. We write $\norm{\bR}$ for the spectral norm of $\bR$, that is the square root of the largest eigenvalue of $\bR\bR^{\top}$. 

For a sample $(\by_1,\ldots, \by_n)=\bSigma^{1/2} \bX$ from the population $\by$ with $\bX=(x_{ij})_{i=1,\ldots,p; j=1,\ldots,n}$, the (non-centered) sample correlation matrix $\hat \bR$ is given by
\begin{equation*}
\hat \bR= \diag(\bS)^{-1/2} \,\bS \, \diag(\bS)^{-1/2}\,,
\end{equation*}
where $\bS= (1/n) \bSigma^{1/2} \bX \bX^\top \bSigma^{1/2}$ is the (non-centered) sample covariance matrix. We are interested in the asymptotic fluctuations of the logarithmic determinant of $\hat \bR$ as $p$ and $n$ tend to infinity simultaneously. 

{\it Throughout the paper, we will assume that the dimension $p$ is a function of the sample size $n$, i.e., $p = p_n$, and that $p/n\to\gamma\in(0, 1]$. All limits are for $\nto$, unless explicitly stated otherwise.} Strictly speaking, the sample correlation and covariance matrices as well as their population counterparts depend on $n$, that is $\hat \bR=\hat \bR_n, \bS=\bS_n,  \bR= \bR_n, \bSigma= \bSigma_n$. For simplicity, we suppress the dependence on $n$ in our notation. We write $\bI$ for the identity matrix if the dimension is clear from the context, and $\cid$ denotes convergence in distribution. \par
\medskip

%\subsection{Main results}

The following CLT is our first main result.

\begin{theorem}[Logarithmic law in the non-centered case]\label{thm:main}
  Assume that $x_{ij}$ are i.i.d.~random variables with mean zero, variance one and finite fourth moment $\mathbbm{E}|x_{11}|^4<\infty$. If the spectral norm of $\bR$ is uniformly bounded and $p/n\to \gamma\in (0, 1]$ with $p< n$, and let
  \begin{align*}
   % \mu_n&=  \log \det (\bR) +\left(p-n+\frac{1}{2}\right)\log\left(1-\frac{p}{n}\right)-p +\frac{p}{n}+ \frac{p}{2n}(\mathbbm{E}|x_{11}|^4-3)\left(C_{\bR^{1/2}}-1 \right) ,\\
		\mu_n&=  \log \det (\bR) +\left(p-n+\frac{1}{2}\right)\log\left(1-\frac{p-1}{n}\right)-(p-1) +\frac{p}{n}+ \frac{p}{2n}(\mathbbm{E}|x_{11}|^4-3)\left(C_{\bR^{1/2}}-1 \right)\,,\\
    \sigma_n^2&=  -2\log\left(1-\frac{p-1}{n}\right)-2\frac{p}{n}+2\frac{p}{n}\tr(\bR-\bI)^2/p\,,
  \end{align*}
  then
  \begin{eqnarray}\label{CLT}
    \frac{\log \det (\hat{\bR})-\mu_n}{\sigma_n}\overset{d}{\longrightarrow}\mathcal{N}(0, 1),\qquad \nto\,,
  \end{eqnarray}
where $C_{\bR^{1/2}}=\frac{1}{p}||\bR^{1/2}\circ\bR^{1/2}||^2_F=\frac{1}{p} \tr\left[\left(\bR^{1/2}\circ\bR^{1/2}\right)^2\right]$, $'\circ'$ denotes the Hadamard product and $\bR^{1/2}$ is the symmetric square root of the matrix $\bR$. 

\end{theorem}
The proof of Theorem~\ref{thm:main} is given in Section~\ref{proofs}. Note that the terms $C_{\bR^{1/2}}$ and $\tr(\bR-\bI)^2/p$ are uniformly bounded in $p$ due to our assumption of boundedness of the largest eigenvalue (spectral norm) of $\bR$. 
The term $\log(1-\tfrac{p-1}{n})$ is kept instead of $\log(1-\tfrac{p}{n})$ in order to incorporate the case $p=n$. In particular, if $p=n$, we deduce from Theorem \ref{thm:main} that \eqref{CLT} holds with
$$\mu_n= \log \det (\bR) -\tfrac{1}{2} \log(n)-n~~~\text{and}~~~\sigma^2_n=2\log n.$$

Noteworthy, if the true mean of the data generating process is known one can already use the above CLT for the purpose of testing but in general one needs to estimate the population mean vector. Thus, one has to consider rather a centered sample correlation matrix
\begin{eqnarray}
  \hat \bR_c= \diag(\bS_c)^{-1/2} \,\bS_c \, \diag(\bS_c)^{-1/2}\,,
\end{eqnarray}
where $\bS_c$ is the centered (by the sample mean) sample covariance matrix given by
\begin{eqnarray}
\bS_c= \frac{1}{n-1} \bSigma^{1/2} (\bX-\bar{\bx}\bi^\top)(\bX-\bar{\bx}\bi^\top)^\top \bSigma^{1/2}~~\text{with $\bar{\bx}=\bX\bi/n$ the sample mean}
\end{eqnarray}
and $\bi=(1,\ldots,1)^{\top}$ denotes the $n$-dimensional vector of ones. In our second main result, we provide a CLT for the logarithmnic determinant of the centered sample correlation matrix.
\begin{theorem}[Logarithmic law in the centered case]\label{sp}
 Assume the conditions of Theorem \ref{thm:main} and $p<n$.  If $p/n\to1$, we assume $p/n=1+O(n^{-1/12})$. Then we have for the logarithmic determinant of the centered sample correlation matrix $\hat{\bR}_c$, 
\begin{eqnarray}\label{CLTcentered}
    \frac{\log \det (\hat{\bR}_c)-\tilde\mu_n}{\sigma_n}\overset{d}{\longrightarrow}\mathcal{N}(0, 1),\qquad \nto\,,
  \end{eqnarray}
where $\tilde\mu_n=\mu_{n-1}$ and $\mu_n, \sigma_n$ are defined as in Theorem~\ref{thm:main}.
\end{theorem}

The proof of Theorem \ref{sp} is given in Section \ref{proofs:sp}. This result is inline with the {\it substitution principle} derived by \cite{zby2015} for linear spectral statistics of sample covariance matrices in the sense that a substitution of $n$ with $n-1$ in the expression for $\mu_n$ yields $\tilde\mu_n$ up to lower order terms.

\begin{remark}{\rm We compare our result with previous ones given in the literature.

(1) It must be noted that the case $\bR=\bI$ in Theorem \ref{sp} was proven for general linear spectral statistics by \cite{gao2017}, one has to, however, compute complex contour integrals first to see the structure of the CLT. The CLT in Theorem \ref{sp} is proven under milder conditions including arbitrary $\bR$ with bounded norm, the case $p/n\to 1$; and it obeys a closed form. We also point out that no information is required about the limiting spectral distribution of the population correlation matrices $\bR$, which appears in the characterization of the Stieltjes transform of the limiting spectral distribution of $\hat\bR_c$; see for example \cite[Theorem~1]{elkaroui:2009}.

(2) Theorem \ref{sp} generalizes the recent result of \cite{jiang2019} to an arbitrary distribution of the entries $x_{ij}$ and removes the restriction on the smallest eigenvalue of $\bR$. Indeed, assuming that $p/n\to\gamma<1$ one can rewrite the limiting mean $\tilde{\mu}_n$ and variance $\sigma^2_n$ in the following way
  \begin{align*}%\label{mean_var_Jiang}
	\tilde{\mu}_n&=  \log \det (\bR) +\left(p-n+\tfrac{3}{2}\right)\log\left(1-\tfrac{p}{n-1}\right)-\frac{n-2}{n-1}p+ \frac{p}{2n}(\mathbbm{E}|x_{11}|^4-3)\left(C_{\bR^{1/2}}-1 \right)+o(1)\,,\\
   \sigma^2_n&=  -2\log\left(1-\frac{p}{n-1}\right)-2\frac{p}{n-1}+\frac{2}{n-1}\tr(\bR-\bI)^2 +o(1)\,,
  \end{align*}
which coincide with the centering and normalization sequences in Theorem 1 of \cite{jiang2019} up to the term $\frac{1}{2}\frac{p}{n}(\mathbbm{E}|x_{11}|^4-3)\left(C_{\bR^{1/2}}-1 \right)$, which is obviously equal to zero in case $x_{ij}\sim\mathcal{N}(0, 1)$ or, more general, $\mathbbm{E}|x_{11}|^4=3$. Nevertheless, it is not the single case when this term disappears, it happens also if $\bR=\bI$. Indeed, by Jensen's inequality $C_{\bR^{1/2}}\geq 1$ with equality if and only if $\bR=\bI$. Interestingly, in this case the statistic on the left-hand sides in \eqref{CLT} and \eqref{CLTcentered} become independent of the moment of order four and one could expect that the restriction $\mathbbm{E}|x_{11}|^4< \infty$ could be weakened\footnote{The investigation of this observation is continued in a subsequent paper.}. 
}\end{remark}

\begin{figure}[h]
 \centering
    \includegraphics[scale=0.45]{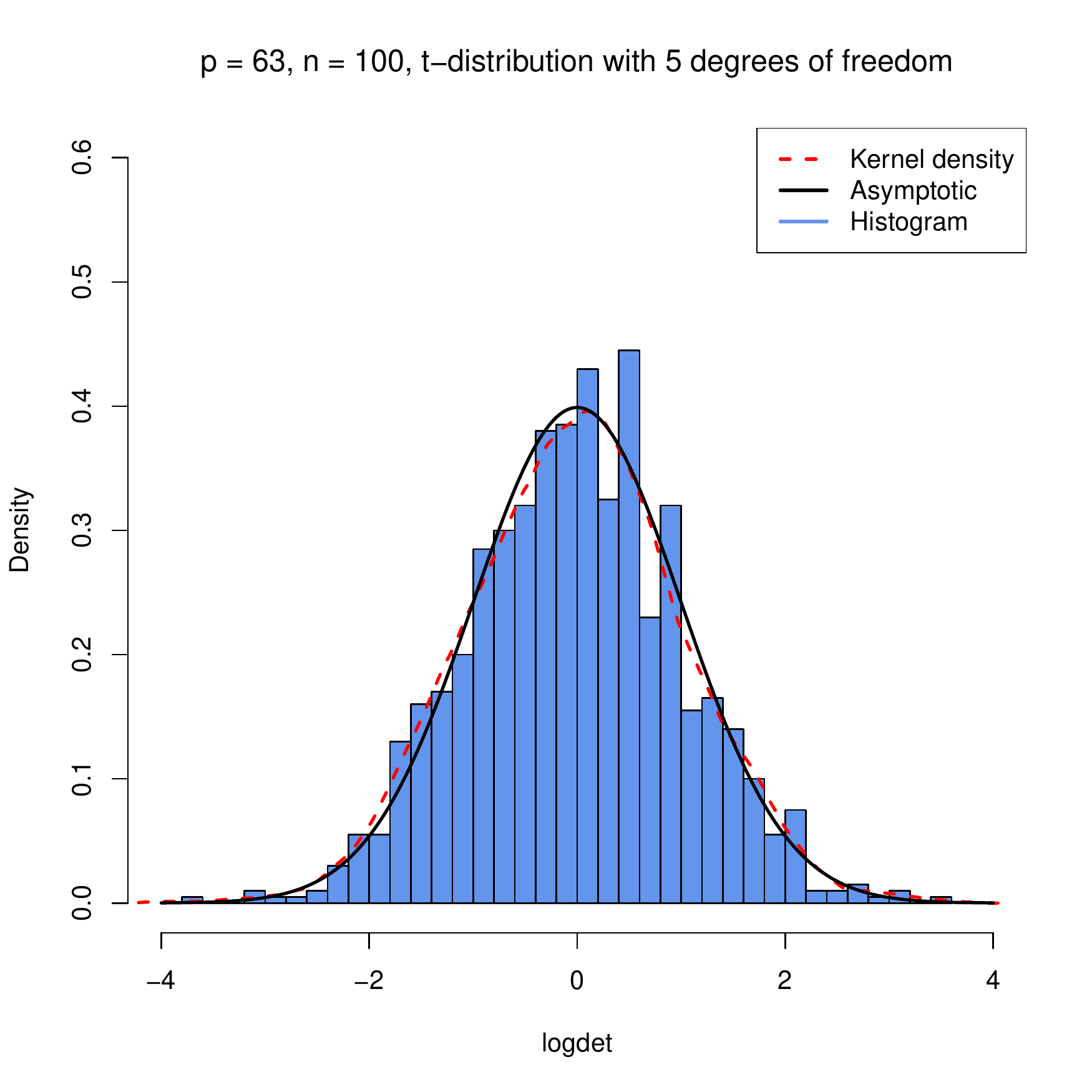} ~~
  \includegraphics[scale=0.45]{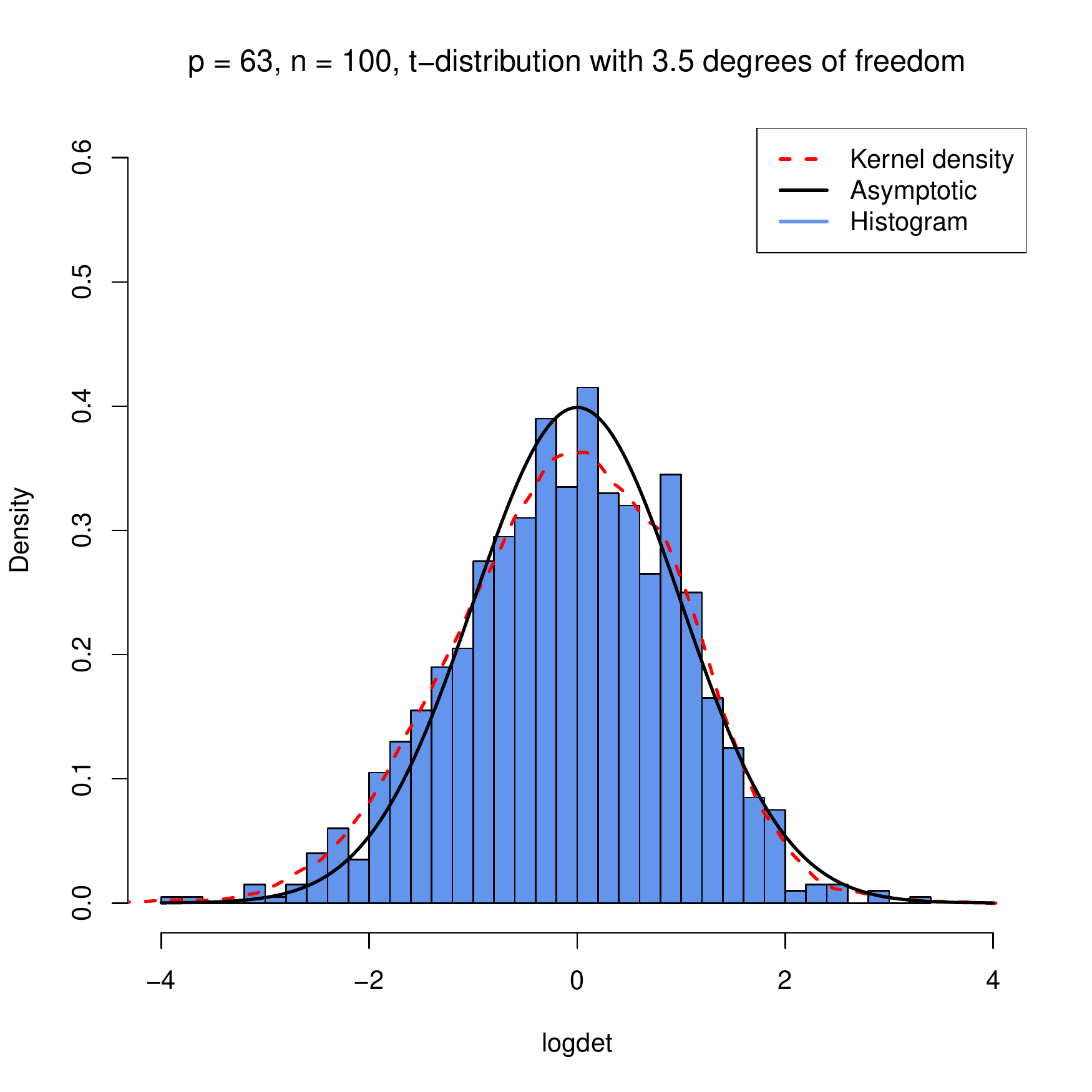}\\
  \includegraphics[scale=0.45]{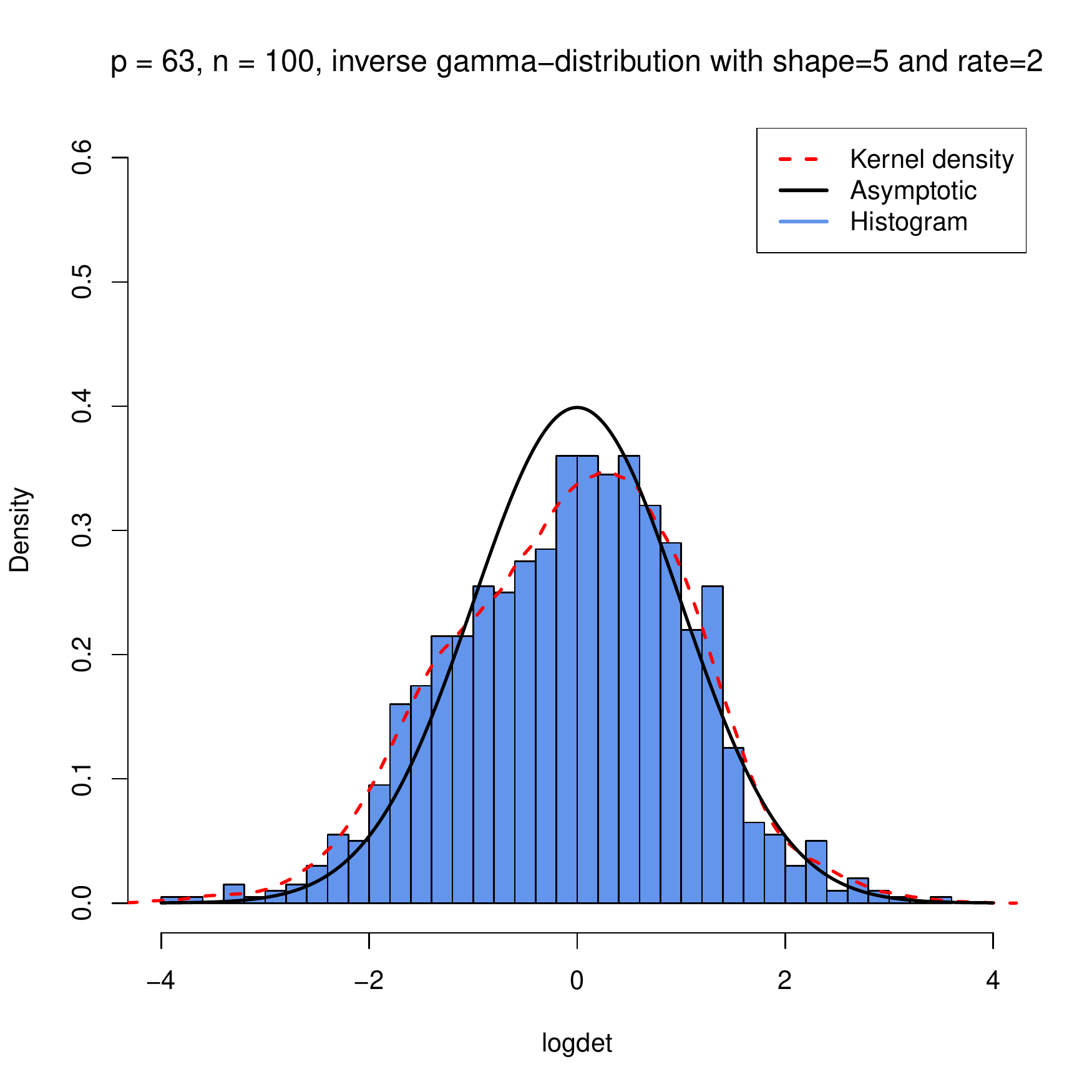}~~
  \includegraphics[scale=0.45]{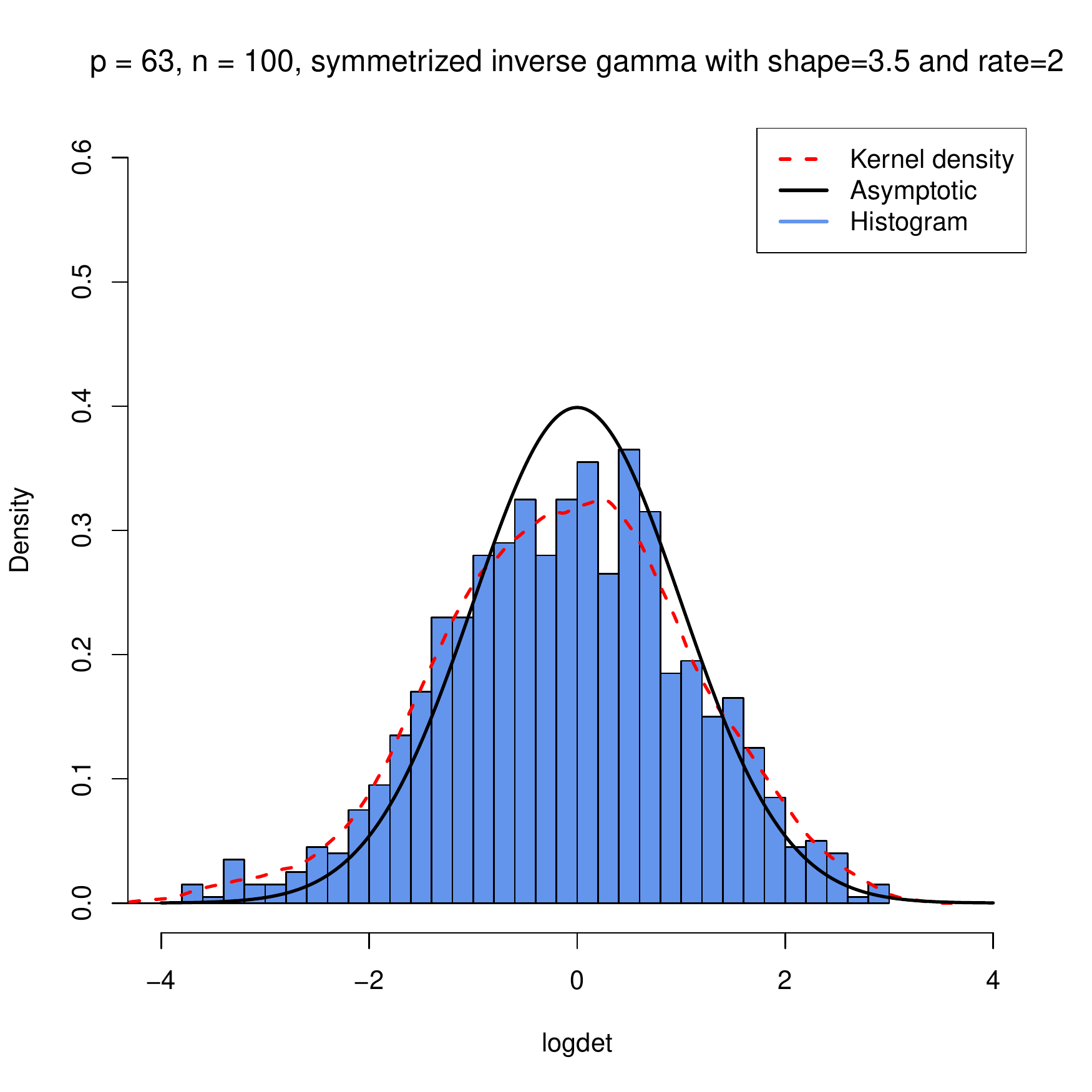}
  \caption{Logarithmic law for $t$ and inverse gamma distribution (left with $\bR=\{0.5^{|i-j|}\}_{i,j=1}^p$ and right with $\bR=\bI$) and $p=63$, $n=100$ with 1000  repetitions.}
  \label{figure1}
\end{figure}

In Figure \ref{figure1}, we have simulated the entries of $\bX$ from a $t$-distribution with different degrees of freedom and from inverse gamma distribution (symmetrized for the case with infinite 4th moment). In case the 4th moment is finite we took the population correlation matrix $\bR=\{0.5^{|i-j|}\}_{i,j=1}^p$, while when the 4th moment is infinite $\bR=\bI$.  We compute $\hat \bR_c$ 1000 times and produce the histogram for $(\log\det(\hat\bR_c)-\tilde{\mu}_n)/\sigma_n$. Then we compare the obtained histogram and kernel density with the standard normal bell curve in order to judge the goodness of fit.

We proceed by discussing Theorem \ref{thm:main} which provides elegant and unified formulas for the asymptotic mean and variance that, in contrast to previous results, avoid heavy computations involving factorials.
The CLT in Theorem \ref{thm:main} is also valid in the case $p/n\to1$ for both $p<n$ and $p=n$, which recently has received particular attention for sample covariance matrices in \cite{ wang2018, nguyen2014, bao2015}. The following remark sheds additional light on this case.

\begin{figure}[h]
  \centering
    \includegraphics[scale=0.45]{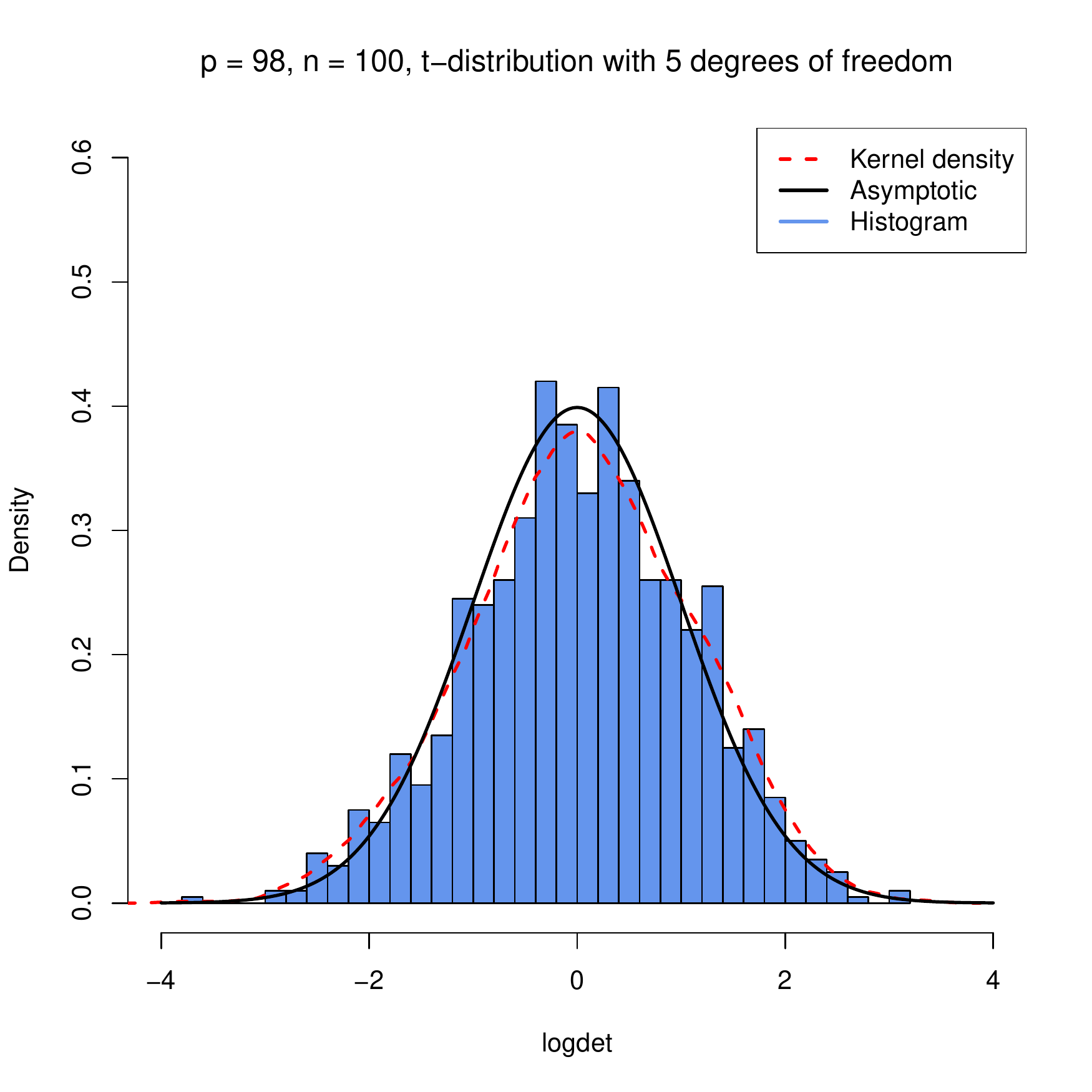} ~~
  \includegraphics[scale=0.45]{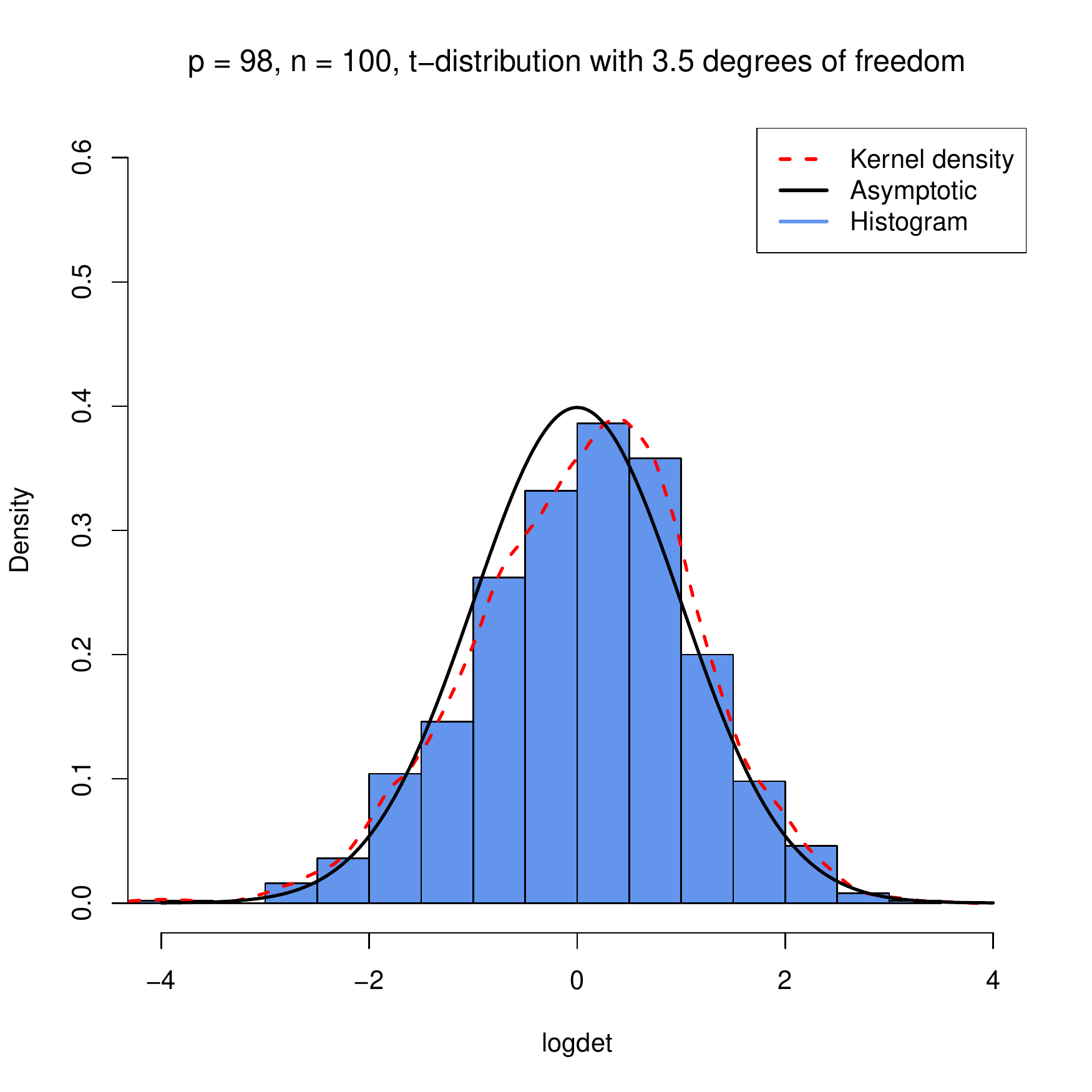}\\
  \includegraphics[scale=0.45]{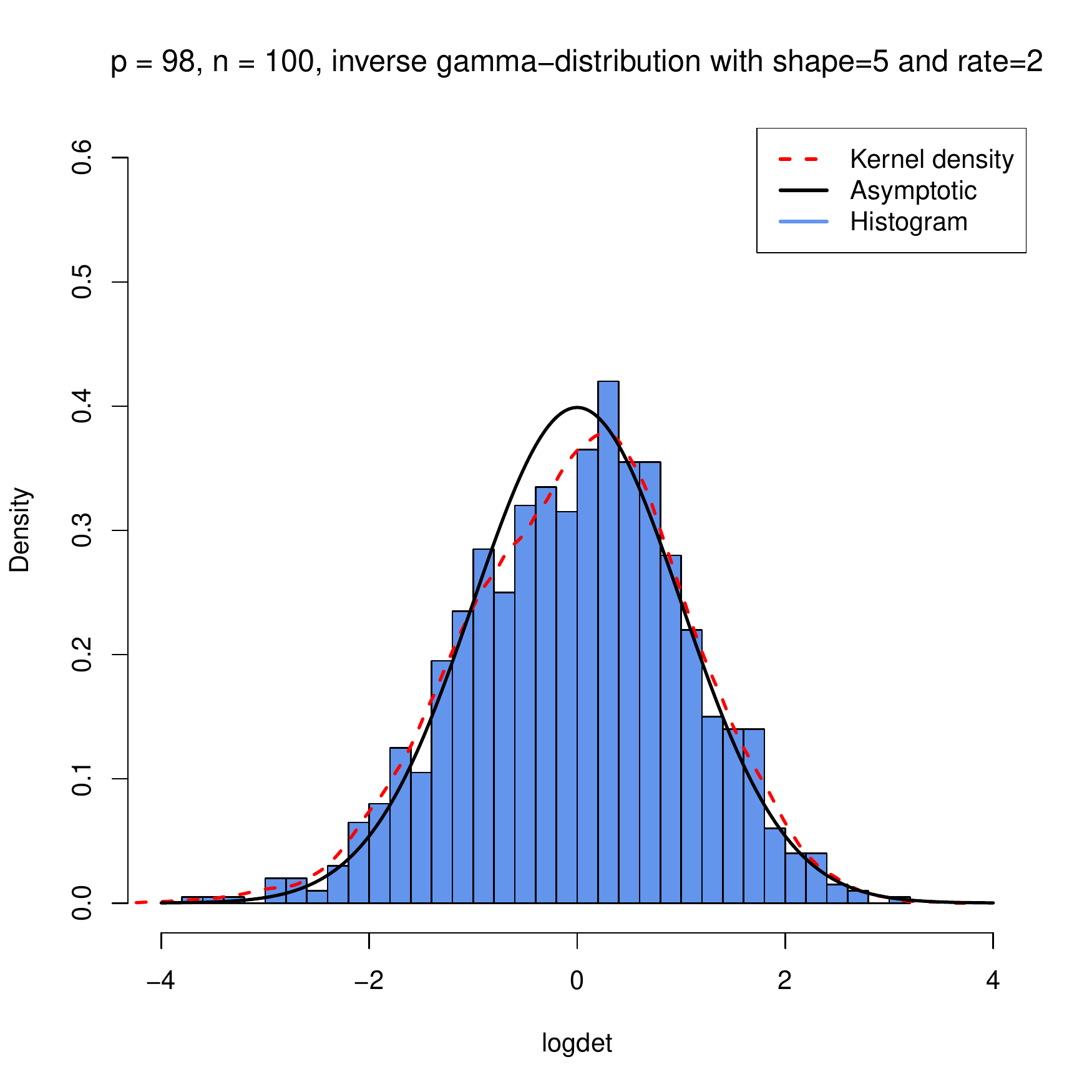}~~
  \includegraphics[scale=0.45]{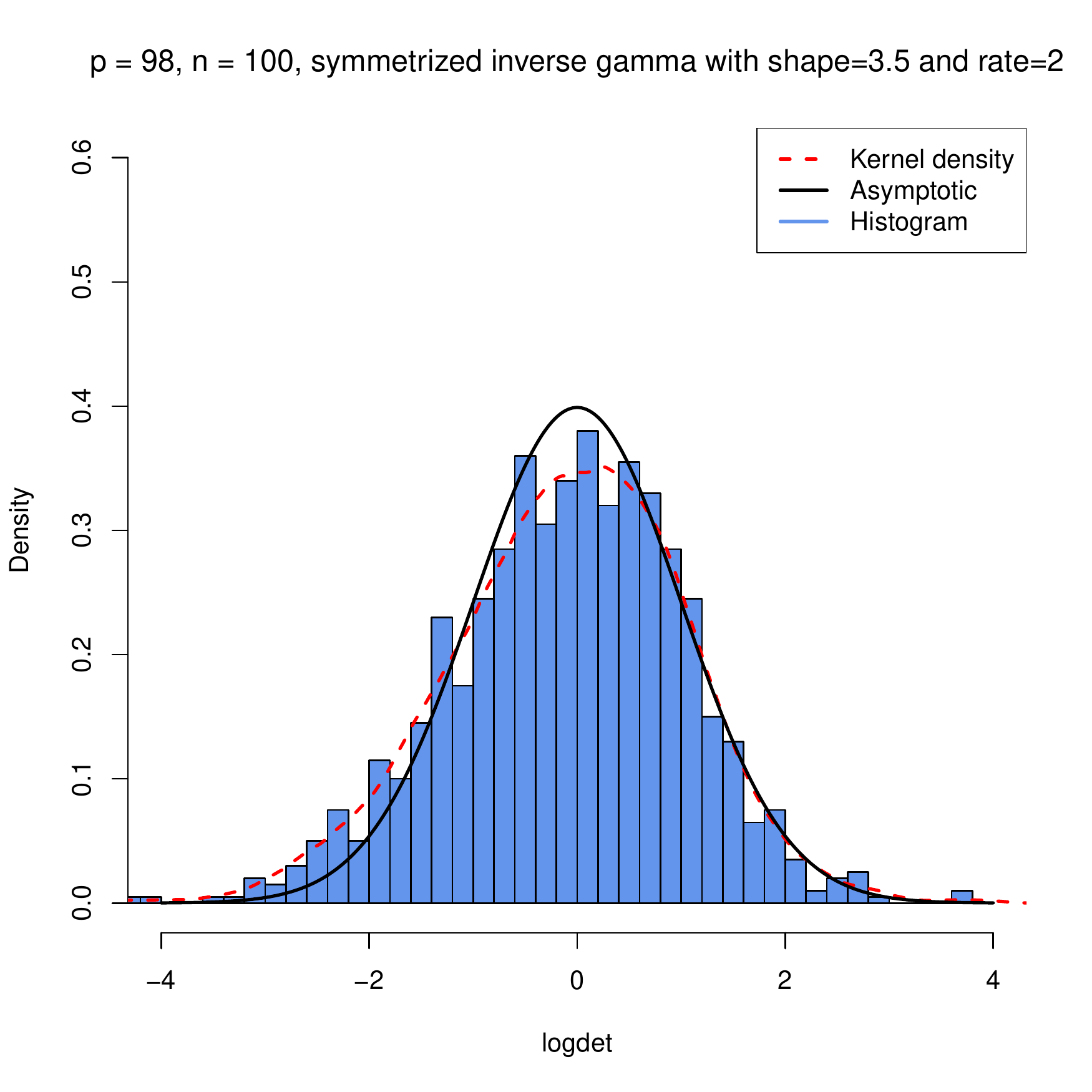}
  \caption{Logarithmic law for $t$ and inverse gamma distribution (left with $\bR=\{0.5^{|i-j|}\}_{i,j=1}^p$ and right with $\bR=\bI$) and $p=98$, $n=100$.}
  \label{figure2}
\end{figure}

\begin{remark}{\rm
(1) We point out that in case $p/n\to 1$ the leading term in variance of the provided CLT tends to infinity. For example, if $n-p$ is constant, $\sigma_n^2$ is of order $\log n$ . Therefore, provided that the largest eigenvalue of $\bR$ is uniformly bounded in $p$, the terms in the asymptotic mean and variance, which are proportional to $p/n$, $(\mathbbm{E}|x_{11}|^4-3)\left(C_{\bR^{1/2}}-1 \right)$ and $\tr(\bR-\bI)^2/p$ will vanish asymptotically. 
 Although one can expect that the convergence to the normal distribution is much slower in this case, this observation has an interesting implication on the test of uncorrelatedness, i.e., $\bR=\bI$. More precisely, it reveals the fact that this test will loose its power asymptotically in case $p/n\to1$ on any alternative hypothesis $\bR\neq\bI$ as long as $\bR$ has no large eigenvalues (spikes). A very similar result was recently found by \cite{bodnar2019} by constructing the test on block diagonality of large covariance matrices.

(2) Interestingly, in case $p/n\to 1$ Theorem \ref{thm:main} also generalizes several results in the literature. First, it is the first result related to the works of \cite{wang2018} and \cite{bao2015}, where the sample covariance matrix was considered. Secondly, one can recover up to the vanishing constants the result of \cite[Theorem 3]{tina2018} taking the result of the CLT for the non-centered sample correlation matrix $\hat{\bR}$, setting $p=n$ and $x_{ij}\sim\mathcal{N}(0, 1)$. The latter we will further use for testing on uniformity on the entries of the large random correlation matrix because the obtained CLT applies also for non-normal data. For illustration we again simulate the entries of the data matrix $\bX$ from $t$ and inverse gamma-distribution 1000 times, similarly as in Figure \ref{figure1}. In Figure \ref{figure2} we plot the kernel densities together with histograms of properly standardized (as to Theorem \ref{sp}) logarithmic determinant of {\it centered} sample correlation matrix $\hat\bR_c$ in case $p=98$ and $n=100$. The asymptotic formula provides still a very convenient fit to the sampled logarithmic determinant. Moreover, the extra assumption $p/n=1+O(n^{-1/12})$ seems to be purely technical.
}\end{remark}

At last, the proof of the main results reveals a very interesting fact about the logarithmic determinant of sample correlation matrix. Indeed, one can show, for example in non-centered case, under generic conditions that the following expansion holds
\begin{eqnarray}\label{asymp_expansion}
 \log \det (\hat{\bR})=\log \det (\bR^{1/2}\tilde{\bS}\bR^{1/2})-\tr(\bR^{1/2}\tilde{\bS}\bR^{1/2}-\bI)+\tfrac{1}{2}\tr(\diag(\bR^{1/2}\tilde{\bS}\bR^{1/2})-\bI)^2+o_{\mathbbm{P}}(1)\,,
\end{eqnarray}
where $\tilde{\bS}=n^{-1} \bX \bX^\top$ and $o_{\mathbbm{P}}(1)$ is a random variable that converges to zero in probability as $\nto$. One may call this {\it asymptotic expansion} of the logarithmic determinant of sample correlation matrix. It shows how the logarithmic determinant of the sample covariance matrix is connected with the latter one for correlation matrix asymptotically. Precisely, going through the proof of Theorem \ref{thm:main} one can show that under the asymptotic regime $p/n\to\gamma\in(0, 1)$, as $\nto$, the first and the second  summands in \eqref{asymp_expansion} are asymptotically jointly normal, whereas the third one is converging to a constant in probability. This result is the key ingredient for further investigations of the centered sample correlation matrix and gives very convenient interpretation of the obtained results.

%\newpage

\section{Applications}\label{sec:3}

In the subsequent sections the derived CLT will be applied for testing the uncorrelatedness of the elements of $\by$ and uniformity of the entries of random correlation matrix.

\subsection{Testing the uncorrelatedness}

Assume that we are interested in testing the hypothesis
\begin{eqnarray}
  H_0: \Corr(\by)=\bI~~\quad \text{vs.}~~\quad H_1: \Corr(\by)\neq \bI\,.
\end{eqnarray}
 In order to provide a proper statistical test we need to construct a test statistic and specify its asymptotic pivotal distribution under $H_0$. In view of Theorem \ref{sp}, a natural test statistic is given by
$T=(\log\det(\hat\bR_c)-\tilde{\mu}_n)/\sigma_n$ which is asymptotically standard normal under $H_0$.
We will reject the null hypothesis $H_0$ in case $T$ is "too large", i.e., larger than the $95\%$ quantile of standard normal distribution.

\begin{figure}[h]

  \vspace{-5cm}

  \centering
   \begin{tabular}{ll}
     \hspace{-2cm}  \includegraphics[scale=0.6]{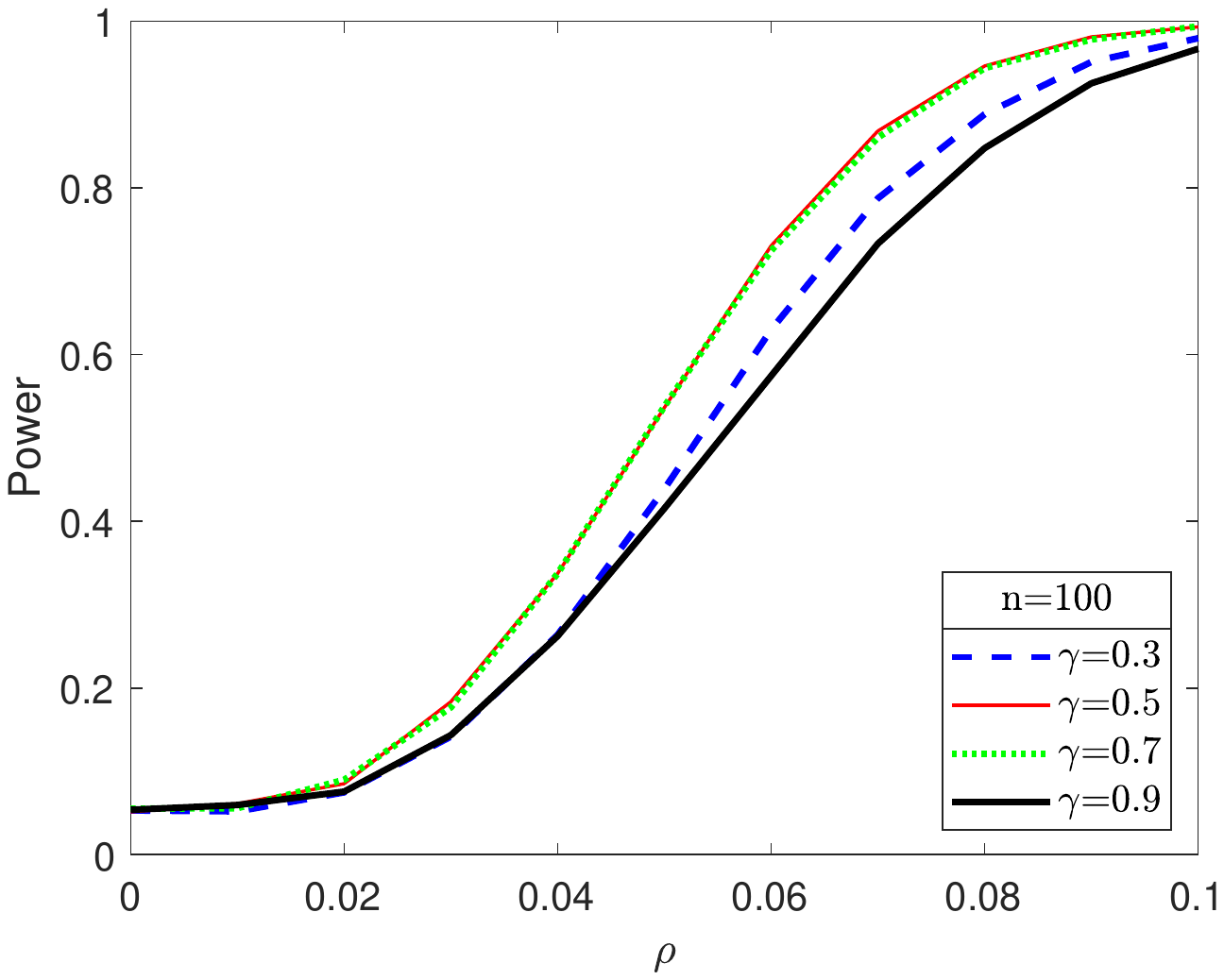} & \hspace{-4cm}\includegraphics[scale=0.6]{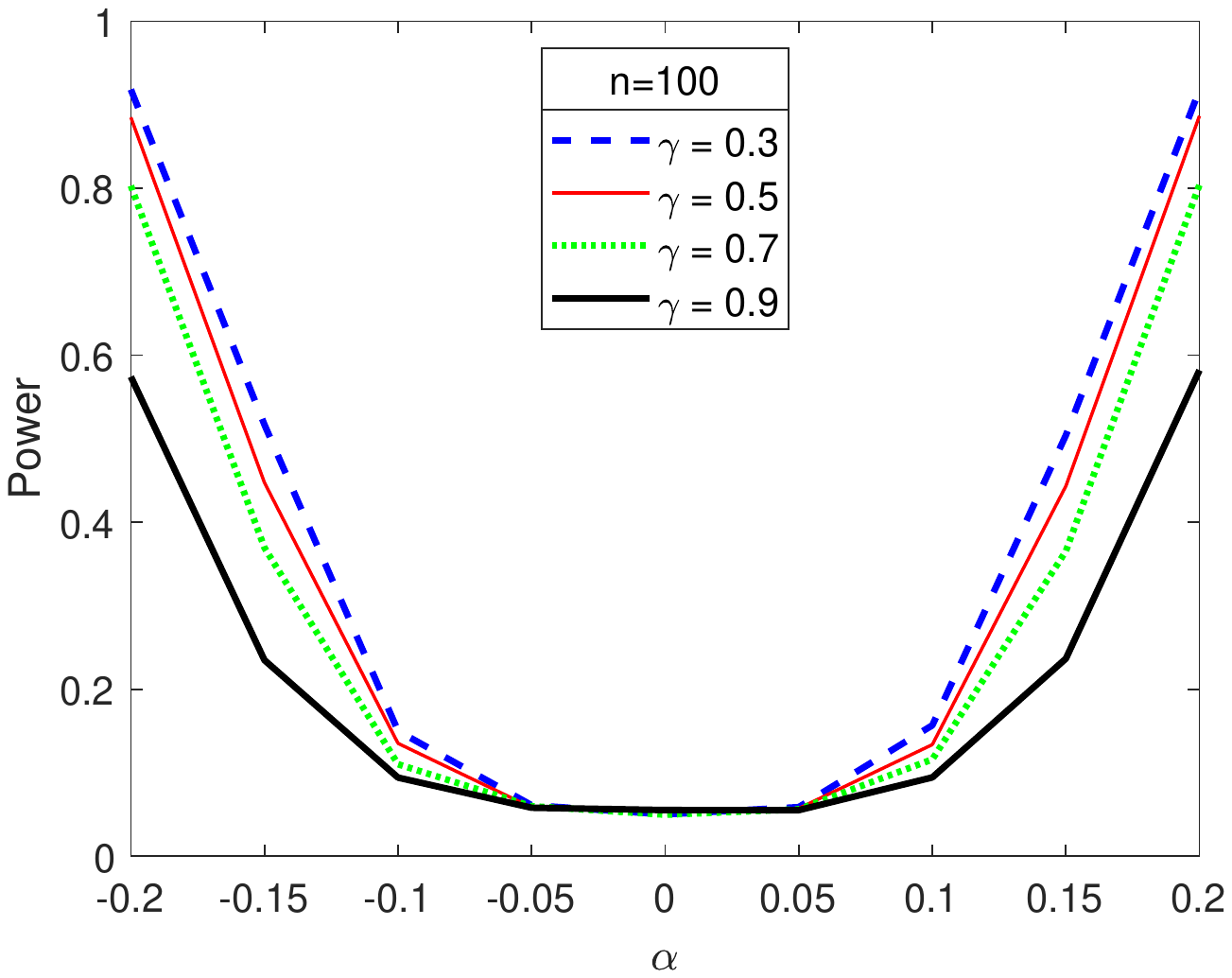}\\[-11.2cm]
      \hspace{-2cm} \includegraphics[scale=0.6]{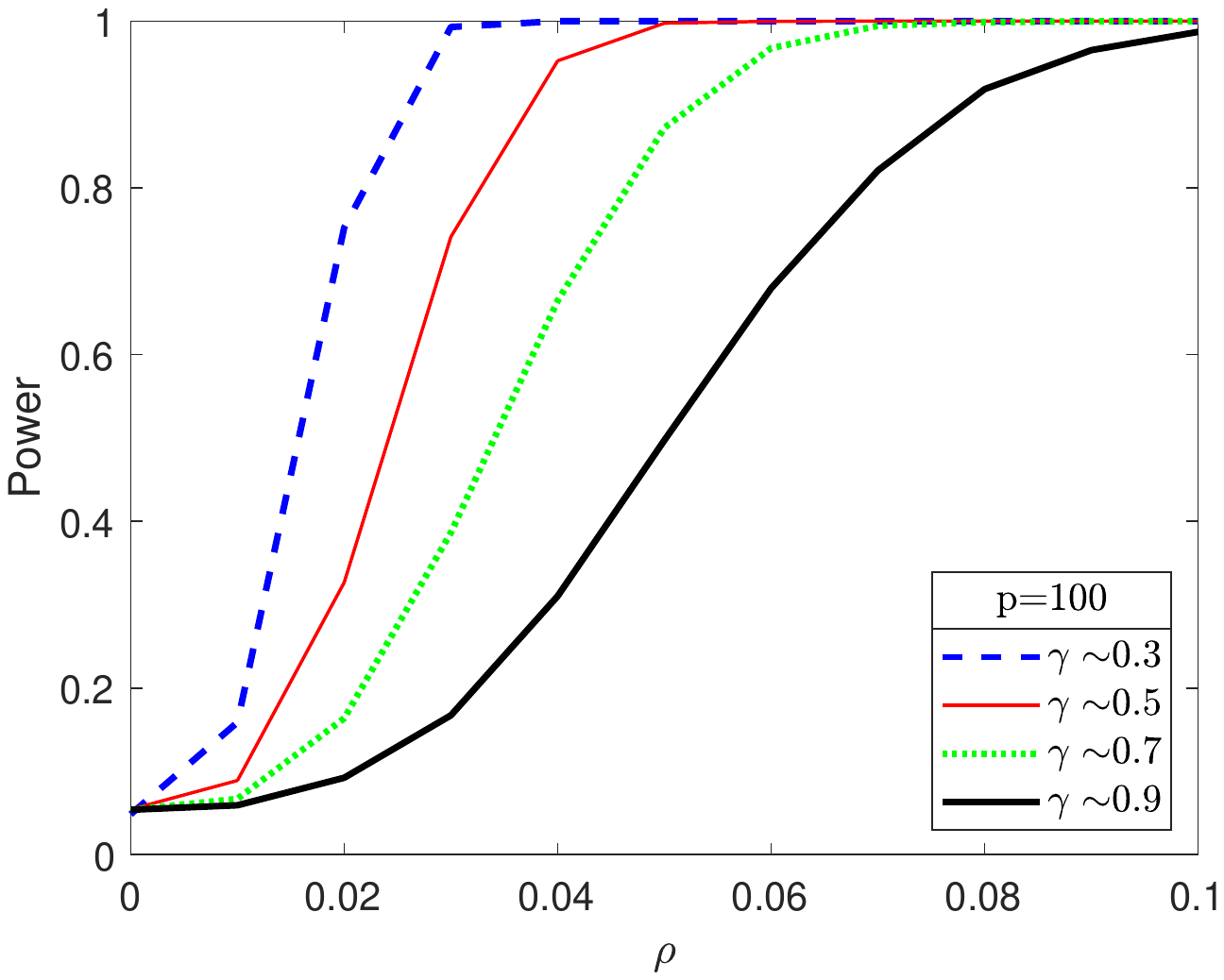}   &    \hspace{-4cm}  \includegraphics[scale=0.6]{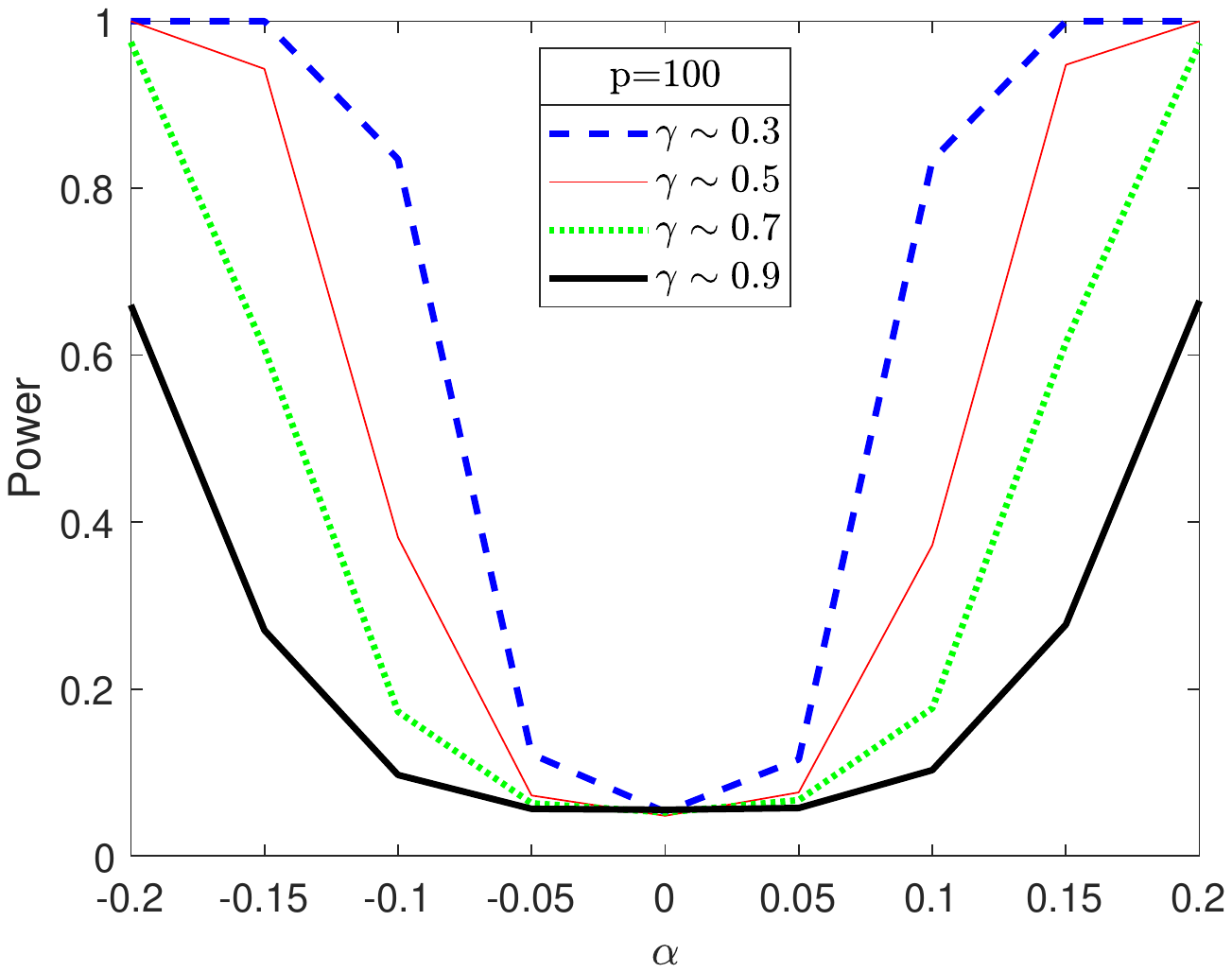}
 \end{tabular}
  \vspace{-5.5cm}
  \caption{Empirical powers of the test on uncorrelatedness for different values of $\gamma=p/n$ for equicorrelation $\bR=(1-\rho)\bI+\rho \bi\bi^\top$ (top) and autocorrelation $\bR=\{ \alpha^{|i-j|}\}_{i,j=1}^p$ (bottom).
     Top: $n$ fixed and $p$ changing, bottom: fixed $p$ and $n$ changing.}
  \label{figure3}
\end{figure}

For that reason we generate the data $(\by_1,\ldots,\by_n)$ with $\by_i=\bSigma^{1/2}\bx_i$, where the components of the noise vector $\bx_i$ are i.i.d. $t$-distributed random variables with $5$ degrees of freedom, mean zero and variance equal to one. The experiment was repeated $N=10^4$ times. Without loss of generality we assume that the population covariance matrix $\bSigma$ is equal to the correlation matrix $\bR$ (standardized data), which is generated in the following manner:

\begin{enumerate}
\item Uncorrelated case: ~~$\bR=\bI$;
\item Autoregressive case:  $\bR=\{ \alpha^{|i-j|}\}_{i,j=1}^p$ for $\alpha\in(-1,1)$;
\item Equicorrelated case:  $\bR=(1-\rho)\bI+\rho \bi\bi^\top$ for $\rho\in(0,1)$.
\end{enumerate}
First, we present the empirical sizes of the uncorrelated case in Table \ref{table1}. The results are plausible even for small values of $(n, p)$ although the size of $5\%$ is a bit overestimated for $n\leq 100$. In general, the larger $p$ the more pronounced is this effect but it vanishes if we increase $p$ and $n$ as expected. Even for $\gamma=0.9$ the test is holding its confidence level quite well, which emphasizes its applicability for moderate finite samples. Next, in Figure \ref{figure3} we present the empirical powers of the test against the equicorrelated case for $\rho\in(0, 0.1)$, where the test shows non-trivial power. On the top part of the Figure \ref{figure3} we consider the equicorrelation case. Left we take $n=100$ and $p$ is varying to get different gamma's from $0.3$ to $0.9$ with a step of $0.2$, while on the right figure we fix $p=100$ and proceed in the same way but with changing $n$. For fixed $n=100$ we observe a tendency towards more power in case of smaller values of $\gamma$ except of the case $\gamma=0.3$, which can be easily explained by too small value of $p=30$, where our asymptotic result seems not to work well. This indicates that $p$ and $n$ should be at least $50$ in this situation to guarantee a reasonable approximation. Indeed, when $p$ is fixed to $100$ we see a natural ranking of the power curves because both $p$ and $n$ are large enough. Interestingly, we have still nontrivial power for $\rho$ as small as $0.08$ even in case when $\gamma$ is near to singularity, i.e., $\gamma=0.9$. A very similar picture is observed on the bottom part of Figure \ref{figure3}, where the power of test against autocorrelation structure was examined. Again, even in the worst case $\gamma=0.9$ the test still provides a reasonable power for $|\alpha|<0.2$. In general, the powers in the autocorrelated case are a bit smaller as in case of equicorrelation, which is not surprising because the correlations are decaying to zero exponentially for the former one.

\vspace{0.5cm}

\begin{table}[h]
    \centering
    \resizebox{0.5\textwidth}{!}{\begin{tabular}[h]{|c||r|r|r|r|}
      \hline
      \backslashbox[4em]{$n$}{$\gamma$} & 0.3 & 0.5 & 0.7 & 0.9\\\hline \hline
    40 &0.0492& 0.0528& 0.0545& 0.0586\\
60 &0.05 &0.0493 &0.0539 &0.0578\\
80 &0.0522 &0.0548 &0.0515 &0.056\\
100 &0.0504 &0.053 &0.053 &0.052\\
500 &0.0501 &0.05 &0.0492 &0.0527\\
      1000 &0.0496 &0.0523 &0.0492 &0.0508 \\   \hline
    \end{tabular}
}
\vspace{0.5cm}

 \caption{Empirical sizes of the test on uncorrelatedness ($\gamma=p/n$ , $N=10000$, $\alpha=0.05$).}
    \label{table1}
  \end{table}

\subsection{Testing the uniformity}

In case of independent normal distributed data, the joint density of the elements of the sample correlation matrix is proportional to its determinant.
More direct (independent of the data generating process) approach of generating a random correlation matrix, such that the density of its entries is proportional to power of the determinant
has been proposed by \cite{joe2006}. In this paper it was shown that any positive definite correlation matrix $\bR$ can be parametrized in terms of appropriately chosen correlations and partial correlations taking independently values from the interval $(-1, 1)$. If these partial correlations are beta distributed (beta distribution transformed to the interval $(-1,1)$) with parameters dependent on the size of the conditioning sets of the partial correlations,  then the joint density of entries of the correlation matrix is proportional to $\det(\bR)^{\eta-1}$, where $\eta>0$.
Each correlation coefficient in such correlation matrix has a $\text{Beta}(\eta-1+p/2, \eta-1+p/2)$ distribution on $(-1,1)$. The uniform joint density is obtained in case $\eta=1$.

Using this fact and the proven CLT we construct a test on the uniformity of the entries of random correlation matrix.
\begin{eqnarray*}
  H_0: \eta=1~ \text{(uniformity)}~~\text{vs.}~H_1: \eta\neq1~\text{(non-uniformity)}\,.
\end{eqnarray*}

\begin{figure}[H]
\vspace{-6.5cm}
  \centering
  \includegraphics[scale=0.7]{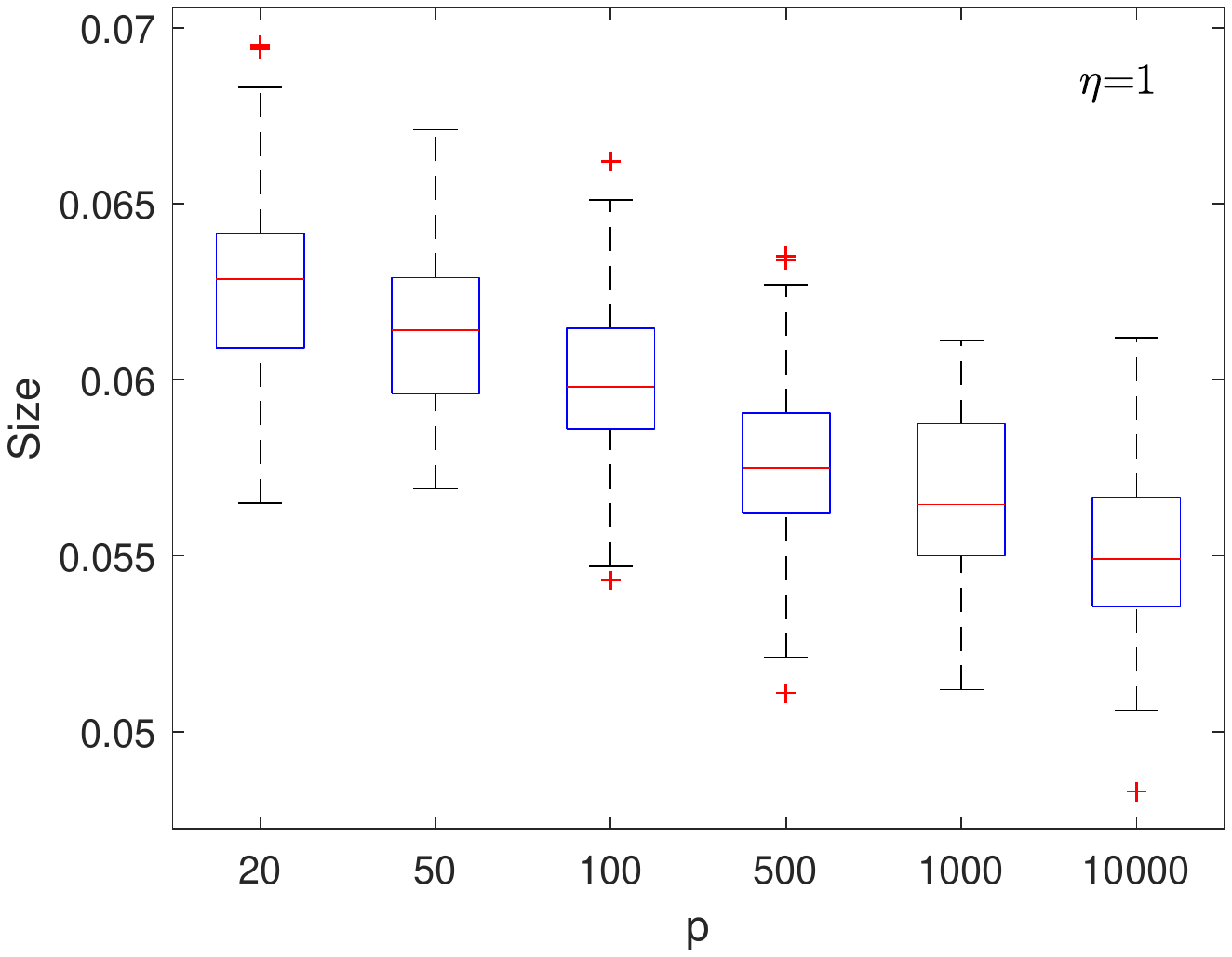}
  \vspace{-6.5cm}
  \caption{Box-plots of empirical sizes of the test on uniformity}
  \label{figure4}
\end{figure}

To test this hypothesis to the level of $5\%$ we generate $N=10000$ random correlation matrices of dimension $p\times p$ using the expansion of the determinant based on partial correlations (due to \cite{joe2006}) for different values of $\eta$ and apply our CLT in Theorem \ref{thm:main} for $n-p=2\eta$. Note that in this situation we do not need to generate data samples, compute correlation matrix and its determinant. We observe only the determinant of $p$-dimensional correlation matrix. Moreover, the formula for the determinant of \cite{joe2006} gives us the possibility to generate really large matrices without loss of efficiency. First, in Figure \ref{figure4} we plot the box-plots of empirical sizes where we observe a slight overestimation of the nominal level of $5\%$ for small dimensions. The convergence to the right size of $0.05$ seems to be quit slow. This is in line with our theoretical finding, where the variance was of order of $\log n$. So, in order to apply our test dimension $p$ must be reasonably large.

\begin{figure}[H]

 \vspace{-5cm}

  \centering
  \begin{tabular}{ll}
     \hspace{-2cm}  \includegraphics[scale=0.6]{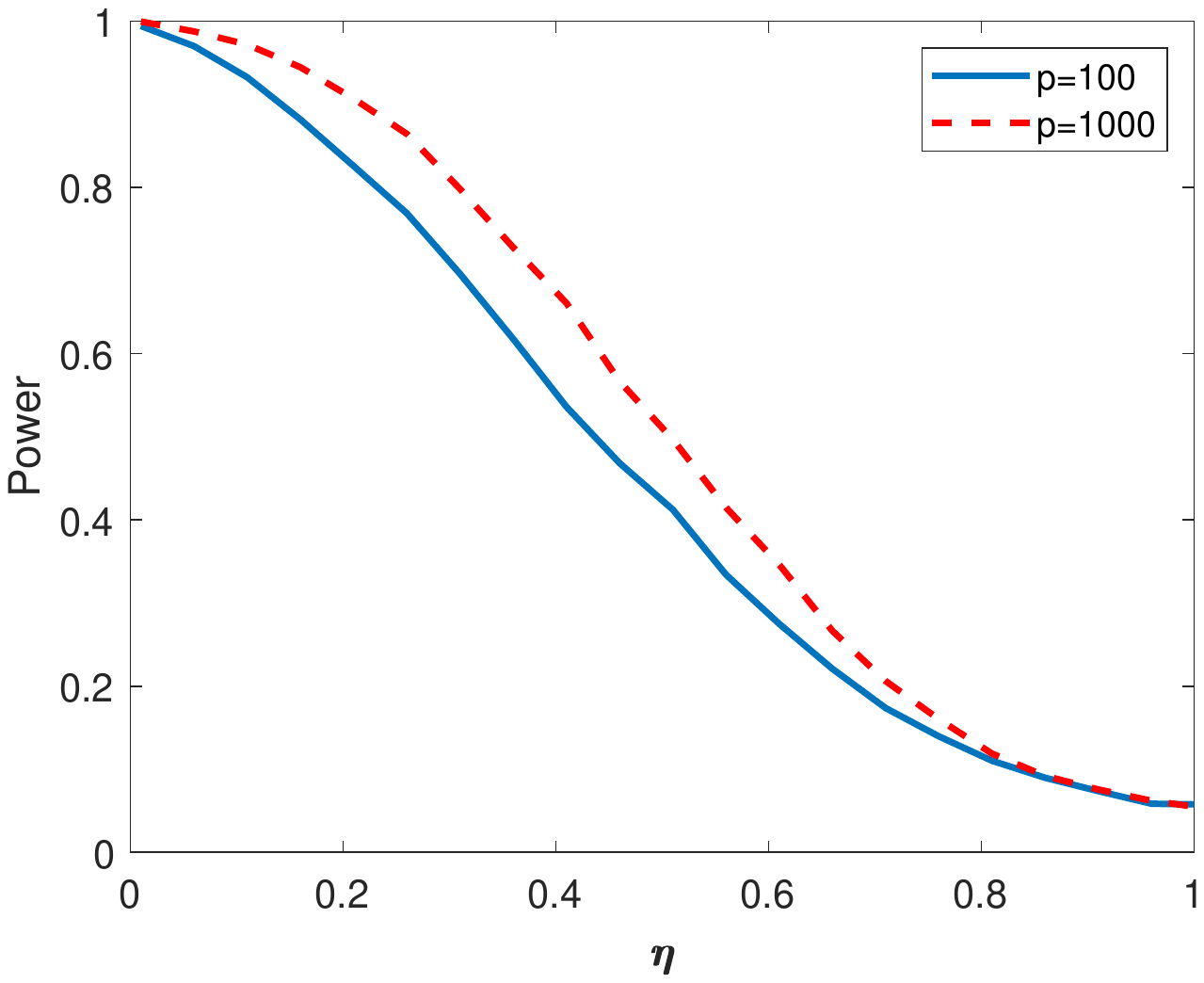} & \hspace{-4cm}\includegraphics[scale=0.6]{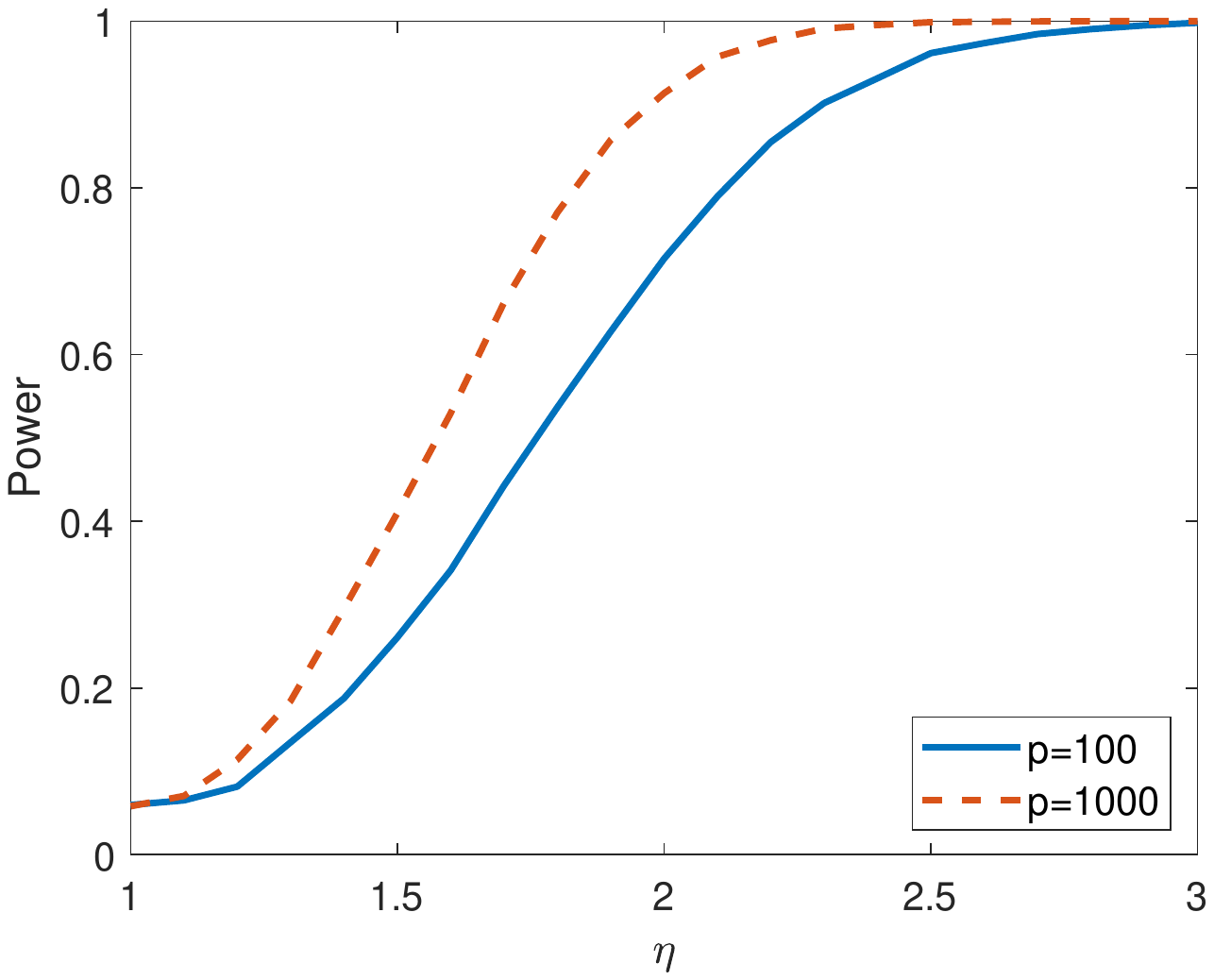}
 \end{tabular}
  \vspace{-5.5cm}
  \caption{Power functions of the test on uniformity for $\eta<1$ (left) and $\eta>1$ (right).}
  \label{figure5}
\end{figure}

\begin{figure}[H]

  \vspace{-6cm}

  \centering
   \begin{tabular}{ll}
     \hspace{-2cm}  \includegraphics[scale=0.6]{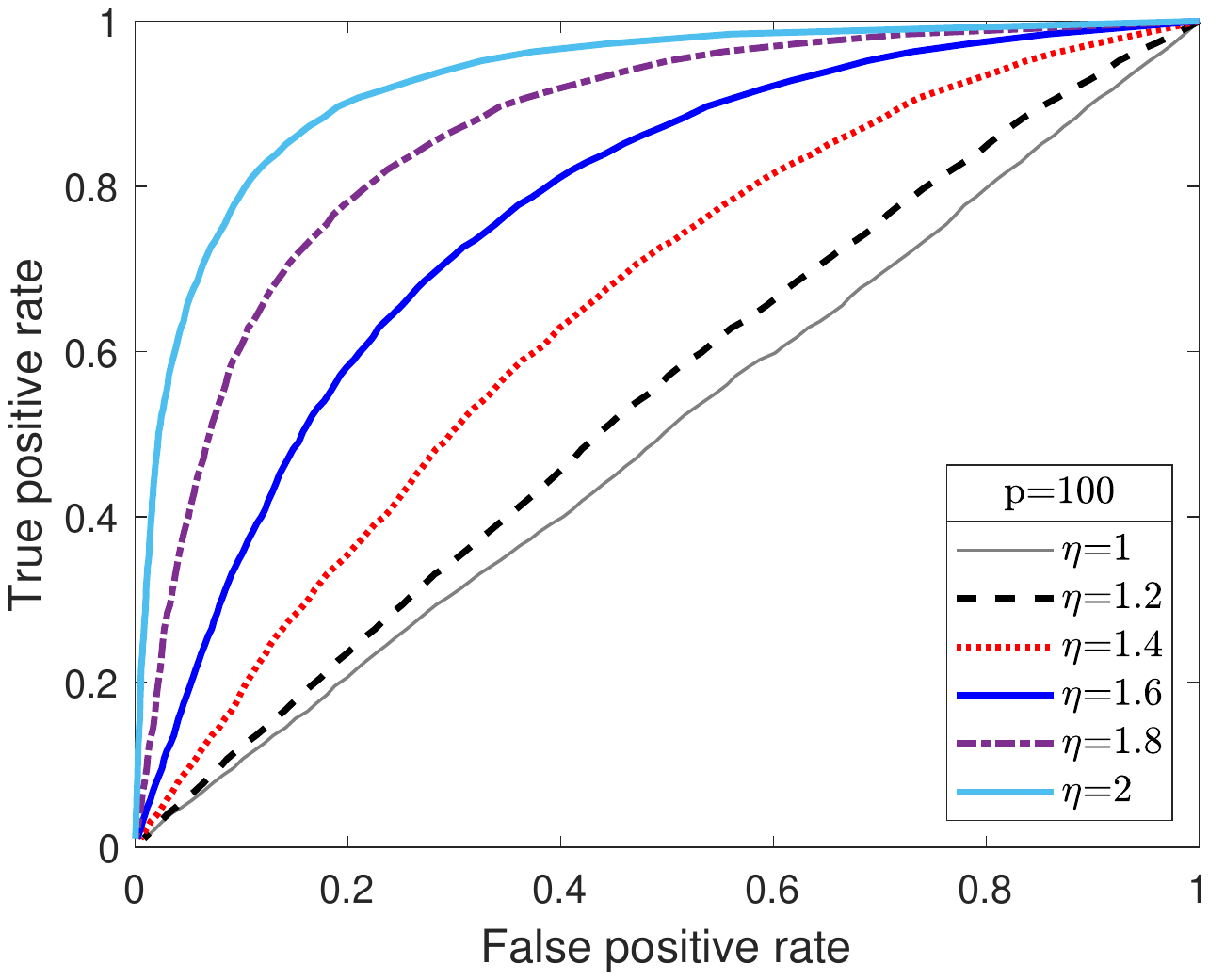} & \hspace{-4cm}\includegraphics[scale=0.6]{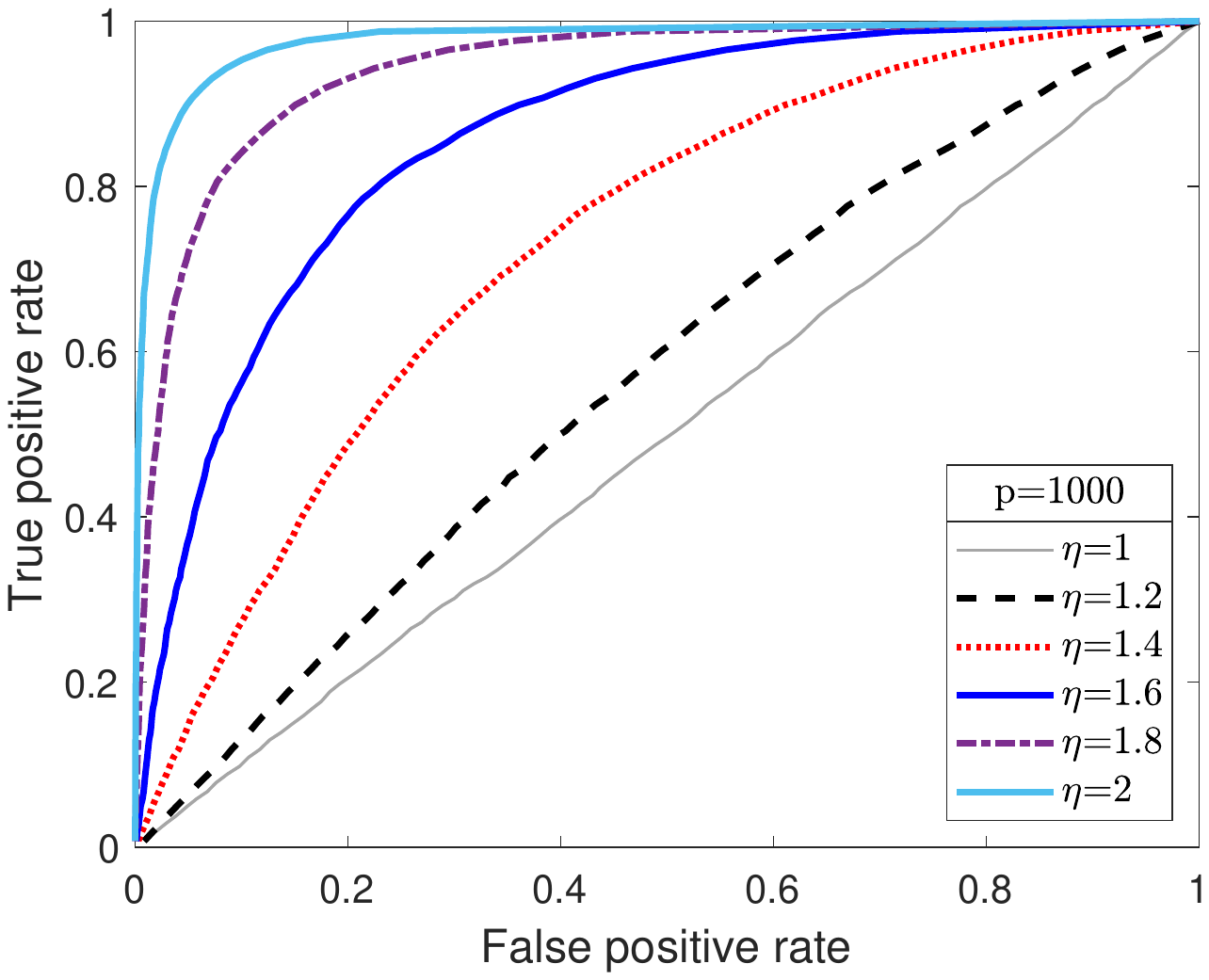}
 \end{tabular}
  \vspace{-5.5cm}
  \caption{ROC curves for the test on uniformity for $p=100$ (left) and $p=1000$ (right) dimensions}.
    %Top: $n$ fixed and $p$ changing, bottom: fixed $p$ and $n$ changing.}
  \label{figure6}
\end{figure}

Next we look on the power of the test against alternative $\eta\neq 1$. The empirical power functions for $p=100$ and $p=1000$ in case of $\eta<1$ and $\eta>1$ are presented in Figure \ref{figure5}. One clearly observes an increasing power when dimension gets larger, e.g., for $\eta>2.5$ ( similarly for $\eta<0.2$) both power curves are close to one. In order to investigate their behavior closer we plot the corresponding ROC (Receiver Operating Characteristic) curves for some fixed values of $\eta>1$. The results are given in Figure \ref{figure6}. The effect of increasing dimension is more pronounced here but the obtained results show an acceptable behavior even for small changes like $\eta=1.4$. This indicates the usefulness of the obtained CLT for the sample correlation matrix in the extreme case $p=n$ for a large class of distributions.


\section*{Conclusions}

In this work we prove the CLT for logarithmic determinant of the large dimensional sample correlation matrix under weak assumptions on data generating process, i.e., we need the existence of the moments of order four. It is also assumed that the dimension of the matrix is proportional to the sample size; the case when they are equal is treated as well. In our work we distinguish between two cases a centered and a non-centered one with slightly different limiting means. Moreover, at the end we apply our obtained results to the testing on uncorrelatedness of high-dimensional random vectors and uniformity of the entries of large random correlation matrix. Simulations suggest that the fourth moment assumption can be further weakened, which opens a new research direction for heavy-tailed sample correlation matrices.

%\newpage

\section{Proofs of the main results}\label{proofs}

{\it We recall that the dimension $p$ is a function of the sample size $n$, i.e., $p = p_n$, and that $p/n\to\gamma\in(0, 1]$. All limits and asymptotic equivalences are for $\nto$, unless explicitly stated otherwise.}

\subsection{Proof of Theorem \ref{thm:main}}

From Section \ref{sec:2} we recall the definition of the population correlation matrix $\bR$  and write $\bX=(\bx_1,\ldots, \bx_n)\in \mathbbm{R}^{p\times n}$.  Denote $\tilde{\bS}=1/n\sum\limits_{i=1}^n\bx_i\bx_i^\top$ and rewrite $\log \det (\hat{\bR})$ in the following way
  \begin{eqnarray}\label{decomp}
    \log \det (\hat{\bR})&=&\log \det (\diag(\bS)^{-1/2}\cdot\bS\cdot\diag(\bS)^{-1/2})=\log \det (\bS)-\log \det (\diag(\bS))\nonumber\\
                            % &=&\text{log det}(\tilde{\bS})+\text{log det}(\bSigma)-\text{log det}(\text{diag}(\bSigma)^{-1/2}\text{diag}(\bS)\text{diag}(\bSigma)^{-1/2})-\text{log det}(\text{diag}(\bSigma))\nonumber\\
                             &=& \log \det (\tilde{\bS}) - \log \det (\diag(\bSigma)^{-1/2}\diag(\bS)\diag(\bSigma)^{-1/2})+\log \det (\bR)\nonumber\\
    &=& \log \det (\bR^{1/2}\tilde{\bS}\bR^{1/2}) - \log \det (\diag(\bR^{1/2}\tilde{\bS}\bR^{1/2}))\,,
  \end{eqnarray}
where the last equality follows from the fact that
\begin{eqnarray*}
\diag(\bSigma)^{-1/2}\diag(\bSigma^{1/2}\tilde{\bS}\bSigma^{1/2})\diag(\bSigma)^{-1/2}
                            &=&\diag\left(\diag(\bSigma)^{-1/2}\bSigma^{1/2}\tilde{\bS}\bSigma^{1/2}\diag(\bSigma)^{-1/2}\right)\\
                            &=&\diag(\bR^{1/2}\tilde{\bS}\bR^{1/2})\,.
\end{eqnarray*}

Now we proceed to the second term in \eqref{decomp}. With the vector $\tilde{\by}_i$ denoting the $i$th row of the matrix $\bR^{1/2}\bX$, it holds
\begin{eqnarray*}
 \log \det (\diag(\bR^{1/2}\tilde{\bS}\bR^{1/2}))=\sum\limits_{i=1}^p\log(\tilde{\by}^\top_{i}\tilde{\by}_i/n)\,,
\end{eqnarray*}
Note that although the elements of one vector $\tilde{\by}_{i}$ are independent zero mean, unit variance variables, the vectors $\tilde{\by}_1,\ldots, \tilde{\by}_p$ depend through the nonzero covariances in the matrix $\bSigma$, i.e., $\text{Cov}(\tilde{\by}_k, \tilde{\by}_l)=\bI_n\bR_{k,l}$ for $k,l=1,\ldots, p$.

Due to \cite{heiny:2021} we have $\max_{i=1,\ldots,p} |\tilde{\by}^\top_{i}\tilde{\by}_i/n-1|\to 0$ almost surely. {\it Recall that all limits are for $\nto$, unless explicitly stated otherwise.}
Therefore, using Taylor expansion of $\log(\tilde{\by}^\top_{i}\tilde{\by}_i/n)=\log(1+(\tilde{\by}^\top_{i}\tilde{\by}_i/n-1))$ around the point zero we get
\begin{eqnarray*}
  \sum\limits_{i=1}^p\log(\tilde{\by}^\top_{i}\tilde{\by}_i/n)&=&\sum\limits_{i=1}^p(\tilde{\by}^\top_{i}\tilde{\by}_i/n-1)-1/2\sum\limits_{i=1}^p(\tilde{\by}^\top_{i}\tilde{\by}_i/n-1)^2+\sum\limits_{i=1}^p\sum\limits_{q=3}^{+\infty}\frac{(-1)^{q+1}}{q}(\tilde{\by}^\top_{i}\tilde{\by}_i/n-1)^q\\
                %                                              &=&\text{tr}(\tilde{\bS}^*-\bI)-\frac{1}{2}\text{tr}(\tilde{\bS}^*-\bI)^2+\sum\limits_{i=1}^pO(|\tilde{\bx}^\top_{i}\tilde{\bx}_i/n-1|^3)\,,\\
                                                              &=&\tr(\bR^{1/2}\tilde{\bS}\bR^{1/2}-\bI)-\frac{1}{2}\tr(\diag(\bR^{1/2}\tilde{\bS}\bR^{1/2})-\bI)^2+\xi_{n,p}\,,
                                                              % \\
       %                                                   &=&-\frac{1}{2}\text{tr}(\text{diag}(\bR^{1/2}\tilde{\bS}\bR^{1/2})^2-\bI)+\xi_{n,p} \,,
\end{eqnarray*}
where $\xi_{n,p}$ is defined by the last equality. Together with \eqref{decomp} this gives

\begin{eqnarray*}
 \log \det (\hat{\bR})=  \log \det (\bR^{1/2}\tilde{\bS}\bR^{1/2})-\tr(\bR^{1/2}\tilde{\bS}\bR^{1/2}-\bI)+\frac{1}{2}\tr(\diag(\bR^{1/2}\tilde{\bS}\bR^{1/2})-\bI)^2-\xi_{n,p}\,.
\end{eqnarray*}

Next we show that $\xi_{n,p}\overset{\mathbbm{P}}{\rightarrow}0$ as $p/n\to \gamma>0$. Because we assumed that $\mathbbm{E}|x_{11}|^4=\nu_4<\infty$ we can use the same truncation method as \cite{baisil2004} page 559, namely we can choose a positive (arbitralily) slowly decreasing to zero sequence $\delta_n$ such that, e.g.,
\begin{eqnarray}\label{delta_n}
  \delta_n\to0,~~\delta_nn^{1/4}\to\infty,~\delta_n^{-4}\mathbbm{E}(x^4_{11}\mathbbm{1}_{\{|x_{11}|>\delta_n\sqrt{n}\}})\to0\,.
\end{eqnarray}
Thus, we may truncate the variables by $\hat{x}_{11}=x_{11}\mathbbm{1}_{\{|x_{11}|\leq \delta_n\sqrt{n}\}}$ without altering the asymptotic results. Note that on the contrary to the original variables all the moments of the truncated variables exist. Indeed, it holds for fixed $n$
\begin{eqnarray}
  \mathbbm{E}|\hat{x}_{11}|^q\leq \delta^{q-4}_nn^{\frac{q}{2}-2}\mathbbm{E}(x^4_{11})\,.
\end{eqnarray}

Using Lemma \ref{bai_silverstein} and \eqref{delta_n} we get for any $2\leq q\leq b\log(n\nu^{-1}_4\delta_n^4)$ and $i=1,\ldots,p$
\begin{eqnarray*}
  \mathbbm{E}\left|\tilde{\by}_i^\top\tilde{\by}_i/n - 1 \right|^q\leq \nu_4(40b^2)^q \delta^{2q-4}_nn^{-1} \,.
\end{eqnarray*}

Denote $s_n=  b\log(n\nu^{-1}_4\delta_n^4)$ and obtain for $\xi_{n,p}$ the following bound
\begin{eqnarray}\label{bound_xi}
  \mathbbm{E}|\xi_{n,p}|&\leq& \sum\limits_{i=1}^p\sum\limits_{q=3}^{+\infty}\frac{\mathbbm{E}\left|\tilde{\by}_i^\top\tilde{\by}_i/n - 1  \right|^q}{q}\leq\frac{1}{3}\lim\limits_{n\to\infty}\sum\limits_{q=3}^{s_n}\sum\limits_{i=1}^p\mathbbm{E}\left|\tilde{\by}_i^\top\tilde{\by}_i/n - 1  \right|^q \nonumber\\
                        &\overset{Lemma~\ref{bai_silverstein}}{\leq}&\lim\limits_{n\to\infty}\sum\limits_{q=3}^{s_n}\frac{p}{3n}\nu_4(40b^2)^q \delta^{2q-4}_n\nonumber\\
                        &\leq&  \underbrace{\lim\limits_{n\to\infty}\frac{p}{3n}\nu_4(40b^2)^2}_{=C<\infty}\lim\limits_{n\to\infty}\sum\limits_{m=2}^{2s_n-4}(\sqrt{40}b\delta_n)^m\leq C\sum\limits_{m=2}^{+\infty}(\sqrt{40}b\delta_n)^m\nonumber\\
                        &=& C\left(\sum\limits_{m=0}^{+\infty}(\sqrt{40}b\delta_n)^m-1-\sqrt{40}b\delta_n\right)=O(\delta^2_n)=o(1)\,,
\end{eqnarray}
where one has to fix $b$ such that $1<b<\frac{\delta_n^{-1}}{\sqrt{40}}$ for $n$ large enough. Now using \eqref{bound_xi} we immediately get that $\xi_{n,p}\cip 0$ for $p/n\to \gamma>0$ and, thus, this result is also valid for $\gamma \in (0, 1]$ and $p\leq n$.

As a result, we have got the following asymptotic expansion of the logarithmic determinant
\begin{eqnarray}\label{logdet_as}
 \log \det (\hat{\bR})=  \log \det (\bR^{1/2}\tilde{\bS}\bR^{1/2})-\tr(\bR^{1/2}\tilde{\bS}\bR^{1/2}-\bI)+\tfrac{1}{2}\tr(\diag(\bR^{1/2}\tilde{\bS}\bR^{1/2})-\bI)^2+o_{\mathbbm{P}}(1)\,.
\end{eqnarray}

The asymptotic equality \eqref{logdet_as} is the key ingredient. It obviously holds
\begin{eqnarray*}
  \log \det (\bR^{1/2}\tilde{\bS}\bR^{1/2})=\log \det (\tilde{\bS})+\log \det (\bR)\,.
\end{eqnarray*}
In the sequel, we will show that the first two summands of the expansion \eqref{logdet_as} are asymptotically jointly normal, whereas the last one converges to some constant in probability. { More precisely, in Section \ref{sec:jointnormal} we will prove for $\gamma\in (0,1)$ that
\begin{equation}\label{eq:jointnormal}
   \frac{ \log \det (\tilde{\bS})-
    \tr(\bR^{1/2}\tilde{\bS}\bR^{1/2}-\bI)-\mu_{\text{log}}}{\sigma_n}
  \cid \mathcal{N}(0,1)
\end{equation}
with $\mu_{\text{log}}=p\left(\frac{p/n-1}{p/n}\log\left(1-\frac{p}{n}\right)-1\right)+\tfrac{1}{2}\log\left(1-\frac{p}{n}\right)-\tfrac{1}{2}(\mathbbm{E}(x^4_{11})-3)\frac{p}{n}$ and 
$$ \sigma_n^2=  -2\log\left(1-\frac{p-1}{n}\right)-2\frac{p}{n}+2\frac{p}{n}\tr(\bR-\bI)^2/p\,. $$
In Section \ref{sec:constant} we will prove that
\begin{equation}\label{eq:constant}
\frac{1}{2}\tr(\diag(\bR^{1/2}\tilde{\bS}\bR^{1/2})-\bI)^2 - \Big[ \frac{1}{2}\frac{p}{n}(\mathbbm{E}(x^4_{11})-3)C_{\bR^{1/2}}+\frac{p}{n} \Big] \cip 0\,.
\end{equation}
Finally, a Taylor expansion yields for $\gamma\in (0,1)$, as $\nto$, 
\begin{equation}\label{eq:seww}
\Big(p-n+\frac{1}{2}\Big)\left[\log\Big(1-\frac{p}{n}\Big)-\log\Big(1-\frac{p-1}{n}\Big) \right] \to -1\,.
\end{equation}
In view of \eqref{logdet_as}, equations \eqref{eq:jointnormal}, \eqref{eq:constant} and \eqref{eq:seww} imply Theorem \ref{thm:main} for $\gamma\in (0,1)$. The case $p/n\to 1$ will be discussed in Section~\ref{sec:pn1}.}

\subsection{Proof of \eqref{eq:jointnormal}}\label{sec:jointnormal}
{
Our goal is to find the asymptotic distribution of $(\log \det (\tilde{\bS}), \tr(\bR^{1/2}\tilde{\bS}\bR^{1/2}-\bI))^\top$ in the case $p/n\to \gamma\in (0, 1)$.

First, consider the term $\log \det (\tilde{\bS})$. 
Now we can apply the CLT proved in \cite{baisil2004} for the sample covariance matrix of noise, namely $\tilde{\bS}$,  with the function $\log(x)$; see also \cite[p.~37]{yao:zheng:bai:2015};  and get for $p/n\to \gamma<1$,
\begin{eqnarray}\label{eq:sdt}
  &&\frac{\log \det (\tilde{\bS}) - \mu_{\text{log}}}{\omega_1} \cid \mathcal{N}(0,1 )
\end{eqnarray}
with $\mu_{\text{log}}$ as in \eqref{eq:jointnormal} and $\omega_1^2=-2\log(1-(p-1)/n)+(\mathbbm{E}(x^4_{11})-3)p/n$. 
Here, we prefer to write $\log(1-(p-1)/n)$ instead of the asymptotically equivalent $\log(1-p/n)$ in order to obtain a unified formula that incorporates the case $p=n$.

For the second term $\tr(\bR^{1/2}\tilde{\bS}\bR^{1/2}-\bI)=\tr[(\tilde{\bS}-\bI)\bR]$ it is easy to check that
\begin{eqnarray}
  \mathbbm{E}\left(\tr(\bR^{1/2}\tilde{\bS}\bR^{1/2}-\bI) \right)&=& \tr (\bR- \bI)=0, \label{eq:exp2nd}\\
  \Var(\tr(\bR^{1/2}\tilde{\bS}\bR^{1/2}-\bI))&=&\mathbbm{E}\left(\tr(\bR^{1/2}\tilde{\bS}\bR^{1/2}-\bI)\right)^2= \frac{1}{n^2}\sum\limits_{i=1}^n \mathbbm{E}\left(\bx^\top_i\bR\bx_i-\tr(\bR)\right)^2\nonumber\\
  &=& \frac{1}{n^2}\sum\limits_{i=1}^n \Var\left(\bx^\top_i\bR\bx_i\right)\nonumber\\
                                                                       &=&  \frac{1}{n^2}\sum\limits_{i=1}^n\left[2\tr\left(\bR^2\right)+(\mathbbm{E}(x^4_{11})-3)\tr(\bR \circ \bR)\right]\nonumber\\
                                                                       %&=& \sum\limits_{i=1}^n\left[2\tr(\bR^2)/n^2+(\mathbbm{E}(x^4_{11})-3)p/n^2\right]\nonumber\\
                                                                       &=& 2\tr(\bR^2)/n-2p/n+(\mathbbm{E}(x^4_{11})-1)p/n \nonumber\\[0.2cm]
  &=&2\tr(\bR-\bI)^2/n+(\mathbbm{E}(x^4_{11})-1)p/n \,,\label{eq:var2nd}
\end{eqnarray}
where for the calculation of the variance we have used Lemma \ref{lem:quf}.

Next, we calculate the covariance.
 Since $\diag(\bR)=\bI$ it holds
\begin{eqnarray}\label{covariance_w12a}
  &&\Cov(\log \det (\tilde{\bS}) ,\tr(\bR^{1/2}\tilde{\bS}\bR^{1/2}-\bI))
=\mathbbm{E}\left(\log \det (\tilde{\bS})\cdot\tr[(\tilde{\bS}-\bI)\bR] \right)\nonumber\\
&& \quad=\mathbbm{E}\left(\log \det (\tilde{\bS})\cdot\tr[(\tilde{\bS}-\bI)] \right)+\mathbbm{E}\left(\log \det (\tilde{\bS})\cdot\tr[\tilde{\bS}(\bR-\bI)] \right)\,.
\end{eqnarray}
Due to \cite{wang2013} we know that as $n\to\infty$ and $p/n\to \gamma\in(0, 1)$ we have
\begin{eqnarray}\label{eq:wang}
  \mathbbm{E}\left(\log \det (\tilde{\bS})\cdot\tr[(\tilde{\bS}-\bI)] \right)= (\mathbbm{E}(x^4_{11})-1)\tfrac{p}{n}+o(1)\,.
\end{eqnarray}

Concerning the second term in \eqref{covariance_w12a}, it is shown in Section~\ref{sec:proof11} that
\begin{equation}\label{eq:prove1}
\lim_{n \to \infty}
\mathbbm{E}\left(\log \det (\tilde{\bS})\cdot\tr[\tilde{\bS}(\bR-\bI)] \right) =0\,.
\end{equation}
Altogether, we obtain from \eqref{covariance_w12a}, \eqref{eq:wang} and \eqref{eq:prove1} that
\begin{equation}\label{eq:cov}
  \Cov(\log \det (\tilde{\bS}) ,\tr(\bR^{1/2}\tilde{\bS}\bR^{1/2}-\bI)) = (\mathbbm{E}(x^4_{11})-1)\tfrac{p}{n}+o(1)
\end{equation}
independently of the structure of population correlation matrix $\bR$.

Now using the result of \cite{najim2016} we get for $\gamma<1$,
\begin{equation}\label{najim_yao}
{\rm d}_{\text{LP}}\left(\mathcal{L} (\log \det (\tilde{\bS})-\mu_{\text{log}}, \tr(\bR^{1/2}\tilde{\bS}\bR^{1/2}-\bI))^\top\,,\, \mathcal{L}(N_n)\right) \to 0\,, \qquad \nto\,,
\end{equation}
where ${\rm d}_{\text{LP}}$ denotes the Levy-Prokhorov distance [we refer to \cite{najim2016} for its definition] and $\mathcal{L}(Y)$ denotes the law of random vector $Y$. The vector $N_n:=(N_{n,1},N_{n,2})^{\top}$ follows a two-dimensional normal distribution with mean vector zero and covariance matrix $\bOmega_n=\left(
  \begin{array}{cc}
    \omega^2_{11,n}  & \omega_{12,n}\\
    \omega_{12,n} & \omega^2_{22,n}
  \end{array}
\right)$. It is important to note that convergence in Levy-Prokhorov distance implies convergence in distribution. 

Next, we determine the entries of $\bOmega_n$. From \eqref{eq:sdt} we get  $\omega_{11,n}^2=-2\log(1-(p-1)/n)+(\mathbbm{E}(x^4_{11})-3)p/n$. Along the lines of the proof of \eqref{eq:var2nd} one can show that $\E\big|\tr(\bR^{1/2}\tilde{\bS}\bR^{1/2}-\bI)\big|^{2+\delta}$ is uniformly bounded in $n$ for some $\delta>0$. Thus, \citet[Theorem 6.2]{DasGupta2008} ensures that also the moments of linear combinations of the components of $(\log \det (\tilde{\bS}), \tr(\bR^{1/2}\tilde{\bS}\bR^{1/2}-\bI))^\top$ converge to the corresponding moments of $N_n$. From \eqref{eq:var2nd} and \eqref{eq:cov}, we deduce that the remaining entries of $\bOmega_n$ can be chosen as
$\omega^2_{22,n}=2\tr(\bR-\bI)^2/n+(\mathbbm{E}(x^4_{11})-1)p/n$ and 
$\omega_{12,n}=(\mathbbm{E}(x^4_{11})-1)\tfrac{p}{n}$ so that

$$\bOmega_n= \left( \begin{array}{cc}
-2\log\left(1-\tfrac{p-1}{n}\right)+(\mathbbm{E}(x^4_{11})-3)\tfrac pn    & (\mathbbm{E}(x^4_{11})-1)\tfrac{p}{n}\\
  (\mathbbm{E}(x^4_{11})-1)\tfrac{p}{n}  & 2\tr(\bR-\bI)^2/n+(\mathbbm{E}(x^4_{11})-1)\tfrac pn \end{array}\right).
$$
Let $\sigma_n$ as in \eqref{eq:jointnormal} and observe that 
\begin{equation*}
   \frac{ N_{n,1}-N_{n,2}}{\sigma_n} \sim \mathcal{N}(0,1)\,.
\end{equation*}
In conjunction with \eqref{najim_yao} this implies \eqref{eq:jointnormal}.
}

\subsubsection{Proof of \eqref{eq:prove1}}\label{sec:proof11}
{
Using $\tilde{\bS}=1/n\sum\limits_{i=1}^n\bx_i\bx_i^\top$ we get
\begin{eqnarray*}
  \mathbbm{E}\left(\log \det (\tilde{\bS})\cdot\tr[\tilde{\bS}(\bR-\bI)] \right)= \frac{1}{n}\sum\limits_{i=1}^n\mathbbm{E}\left(\log \det (\tilde{\bS})\cdot\bx_i^\top(\bR-\bI)\bx_i\right)\,.
\end{eqnarray*}
By the matrix determinant lemma we may write
\begin{eqnarray*}
  \log \det (\tilde{\bS})= \log \det \Big(\underbrace{\tilde{\bS}-\frac{1}{n}\bx_i\bx^\top_i}_{=:\tilde{\bS}_{(i)}}+\frac{1}{n}\bx_i\bx_i^\top \Big)=\log\Big(1+\frac{1}{n}\bx_i^\top\tilde{\bS}^{-1}_{(i)}\bx_i\Big)+ \log \det (\tilde{\bS}_{(i)})\,,
\end{eqnarray*}
which, because of the independence of $\bx_i$ and $\tilde{\bS}_{(i)}$ and $\mathbbm{E}(\bx_i^\top(\bR-\bI)\bx_i)=\tr(\bR-\bI)=0$, implies
\begin{eqnarray}\label{}
   \mathbbm{E}\left(\log \det (\tilde{\bS})\cdot\tr[\tilde{\bS}(\bR-\bI)] \right)= \frac{1}{n}\sum\limits_{i=1}^n\mathbbm{E}\left(\log\left(1+\frac{1}{n}\bx_i^\top\tilde{\bS}^{-1}_{(i)}\bx_i\right)\cdot\bx_i^\top(\bR-\bI)\bx_i\right)\,.
\end{eqnarray}
Next we consider the term $\log\left(1+\frac{1}{n}\bx_i^\top\tilde{\bS}^{-1}_{(i)}\bx_i\right)$. It holds
\begin{eqnarray*}
%&& \log\left(1+\frac{1}{n}\bx_i^\top\tilde{\bS}^{-1}_{(i)}\bx_i\right)\\
%&=& \log\left(1+\frac{1}{n}\bx_i^\top\tilde{\bS}^{-1}_{(i)}\bx_i\right)- \log\left(1+\frac{1}{n}\mathbbm{E}(\bx_i^\top\tilde{\bS}^{-1}_{(i)}\bx_i)\right)+\log\left(1+\frac{1}{n}\mathbbm{E}(\bx_i^\top\tilde{\bS}^{-1}_{(i)}\bx_i)\right)\\
\log\left(1+\frac{1}{n}\bx_i^\top\tilde{\bS}^{-1}_{(i)}\bx_i\right)&=& \log\left(\frac{1+\frac{1}{n}\bx_i^\top\tilde{\bS}^{-1}_{(i)}\bx_i}{1+\frac{1}{n}\tr(\tilde{\bS}^{-1}_{(i)})}\right)+\log\left(1+\frac{1}{n}\tr(\tilde{\bS}^{-1}_{(i)})\right)\\
&=& \log\Big(1+\underbrace{\frac{\frac{1}{n}\bx_i^\top\tilde{\bS}^{-1}_{(i)}\bx_i-\frac{1}{n}\tr(\tilde{\bS}^{-1}_{(i)})}{1+\frac{1}{n}\tr(\tilde{\bS}^{-1}_{(i)})}}_{=:Z_i}\Big)+\log\left(1+\frac{1}{n}\tr(\tilde{\bS}^{-1}_{(i)})\right)\,.
%&=& \log\left(1+Z_i\right)+\log\left(1+\frac{1}{n}\tr(\tilde{\bS}^{-1}_{(i)})\right)\,.
\end{eqnarray*}
This implies that
\begin{eqnarray*}
  \mathbbm{E}\left(\log \det (\tilde{\bS})\cdot\tr[\tilde{\bS}(\bR-\bI)] \right)&=&\frac{1}{n}\sum\limits_{i=1}^n\mathbbm{E}\left[ \left(\log\left(1+{Z}_i\right)+\log\left(1+\frac{1}{n}\tr(\tilde{\bS}^{-1}_{(i)})\right)\right)\cdot\bx_i^\top(\bR-\bI)\bx_i\right]\\
&=& \frac{1}{n}\sum\limits_{i=1}^n\mathbbm{E}\left[ \left(\log\left(1+{Z}_i\right)\right)\cdot\bx_i^\top(\bR-\bI)\bx_i\right]\,.
\end{eqnarray*}
Now, using Lemma \ref{bai_silverstein} we can see that
\begin{eqnarray*}
\mathbbm{E}|Z_i|^4\le  \mathbbm{E}\left|\frac{1}{n}\bx_i^\top\tilde{\bS}^{-1}_{(i)}\bx_i-\frac{1}{n}\tr(\tilde{\bS}^{-1}_{(i)})\right|^4\leq C \,\E \big[\norm{\tilde{\bS}^{-1}_{(i)}}^4\big] \, \delta_n^{4}n^{-1}= O(\delta_n^{4}n^{-1})\,,
\end{eqnarray*}
because $\norm{\tilde{\bS}^{-1}_{(i)}} \to (1-\sqrt{\gamma})^{-2}<\infty$. By the union bound and Markov's inequality, we have for any small $\varepsilon>0$,
\begin{eqnarray*}
  \mathbbm{P}\left(\max\limits_{i=1,\ldots,n}|Z_i|>\varepsilon\right)\leq \sum\limits_{i=1}^n  \mathbbm{P}\left(|Z_i|>\varepsilon\right)\le \frac{1}{\varepsilon} \sum\limits_{i=1}^n\mathbbm{E} |Z_i|^4 = O(\delta_n^4)=o(1)\,,
\end{eqnarray*}
which implies that $\max_{i=1,\ldots,n}|Z_i|\cip 0$ as $\nto$.
The latter justifies the Taylor expansion of the logarithm, namely it holds that
\begin{eqnarray*}
  \log\left(1+{Z}_i\right)=Z_i-\frac{Z_i^2}{2}+\sum\limits_{q=3}^{+\infty} \frac{(-1)^{q+1}}{q}Z_i^q\,, \qquad i=1,\ldots,n\,.
\end{eqnarray*}
Thus, we get
\begin{eqnarray*}
 && \frac{1}{n}\sum\limits_{i=1}^n\mathbbm{E}\left[\log\left(1+Z_i\right)\cdot\bx_i^\top(\bR-\bI)\bx_i\right]\\
&=&\frac{1}{n}\sum\limits_{i=1}^n\mathbbm{E}\left[Z_i\cdot\bx_i^\top(\bR-\bI)\bx_i\right]-\frac{1}{2n}\sum\limits_{i=1}^n\mathbbm{E}\left[Z^2_i\cdot\bx_i^\top(\bR-\bI)\bx_i\right]+\theta_n\,,
\end{eqnarray*}
where $\theta_n=\frac{1}{n}\sum_{i=1}^n\E\left[\log\left(1+Z_i\right)\cdot\bx_i^\top(\bR-\bI)\bx_i\right]$.
Regarding the first term we have using Lemma \ref{lem:quf}
\begin{eqnarray*}
  \mathbbm{E}\left[Z_i\cdot\bx_i^\top(\bR-\bI)\bx_i\right]&=&\mathbbm{E}\left(\left.\mathbbm{E}\left[\frac{\frac{1}{n}\bx_i^\top\tilde{\bS}^{-1}_{(i)}\bx_i-\frac{1}{n}\tr(\tilde{\bS}^{-1}_{(i)})}{1+\frac{1}{n}\tr(\tilde{\bS}^{-1}_{(i)})}\cdot\bx_i^\top(\bR-\bI)\bx_i\right|\tilde{\bS}_{(i)}\right]\right)\\
&=& \mathbbm{E}\left(\frac{1}{1+\frac{1}{n}\tr(\tilde{\bS}^{-1}_{(i)})}\left.\mathbbm{E}\left[\left(\frac{1}{n}\bx_i^\top\tilde{\bS}^{-1}_{(i)}\bx_i-\frac{1}{n}\tr(\tilde{\bS}^{-1}_{(i)})\right)\cdot\bx_i^\top(\bR-\bI)\bx_i\right|\tilde{\bS}_{(i)}\right]\right)\\
&=&\mathbbm{E}\left(\frac{\frac{2}{n}\tr(\tilde{\bS}^{-1}_{(i)}(\bR-\bI))}{1+\frac{1}{n}\tr(\tilde{\bS}^{-1}_{(i)})}\right)\approx\frac{2}{1-\gamma}\frac{\tr(\bR-\bI)}{1+\frac{\gamma}{1-\gamma}}n^{-1}+O(n^{-1})=O(n^{-1})\,.
\end{eqnarray*}

The symbol ''$\approx$'' in the last line means asymptotic equivalence in terms of the convergence of absolute difference to zero. This result follows from \citet[Lemma A1]{Bodnar2018}\footnote{It has to be noted that \cite{Bodnar2018} need here $4+\varepsilon$ moments to exist, while we assume $\mathbbm{E}|x_{11}|^4<\infty$. This is due to the fact that they consider the almost sure convergence, while in our situation the convergence in probability is enough. Indeed, it can be easily shown that this extra $\varepsilon$ follows from the Borel-Cantelli lemma.} and convergence is uniform over $i$ because $\tilde{\bS}_{(i)}$ can be safely replaced by $\tilde{\bS}$ without altering the limit. Similarly, for the second term we have again due to Lemma \ref{lem:quf} and \citet[Lemma A1 and (A.18)]{Bodnar2018}
\begin{eqnarray*}
 \mathbbm{E}\left[Z^2_i\cdot\bx_i^\top(\bR-\bI)\bx_i\right]&=&\mathbbm{E}\left(\frac{\frac{4}{n^2}\tr(\tilde{\bS}^{-1}_{(i)})\tr(\tilde{\bS}^{-1}_{(i)}(\bR-\bI))+\frac{8}{n^2}\tr(\tilde{\bS}^{-2}_{(i)}(\bR-\bI))}{(1+\frac{1}{n}\tr(\tilde{\bS}^{-1}_{(i)}))^2}\right)\\
&\approx& \frac{\frac{4\gamma}{(1-\gamma)^2}\tr(\bR-\bI)n^{-1}+\frac{8}{(1-\gamma)^3}\tr(\bR-\bI)n^{-2}}{\left(1+\frac{\gamma}{1-\gamma}\right)^2}  +O(n^{-1})=O(n^{-1})\,.
\end{eqnarray*}

Thus, to finish the proof of \eqref{eq:prove1} it is enough to show that $\theta_n\overset{\mathbbm{P}}{\to}0$. Indeed, similarly to \eqref{bound_xi} we use Cauchy-Schwarz (CS) inequality, the uniform boundedness of $\norm{\bR}$ and Lemma \ref{bai_silverstein} to receive
\begin{eqnarray*}
  |\theta_n|&\leq& \frac{1}{n}\sum\limits_{i=1}^n\sum\limits_{q=3}^{+\infty}\mathbbm{E}(|Z_i|^q\cdot\left|\bx_i^\top(\bR-\bI)\bx_i\right|)\overset{CS}{\leq} \frac{1}{n} \sum\limits_{i=1}^n\sum\limits_{q=3}^{+\infty}\mathbbm{E}^{1/2}(|Z_i|^{2q})\cdot\mathbbm{E}^{1/2}(\left|\bx_i^\top(\bR-\bI)\bx_i\right|^2)\\
&\leq&\frac{1}{\sqrt{n}} \sum\limits_{i=1}^n\sum\limits_{q=3}^{+\infty}\mathbbm{E}^{1/2}(|Z_i|^{2q})\cdot\left(\frac{1}{n}\tr(\bR-\bI)^2\right)^{1/2}\\
&\overset{\text{Lemma \ref{bai_silverstein}}}{\leq}& \frac{C}{\sqrt{n}} \sum\limits_{i=1}^n\sum\limits_{q=3}^{+\infty}\nu_4^{1/2}\left(\frac{\sqrt{40}b}{(1-\sqrt{\gamma})}\right)^q\delta_n^{q-2}n^{-1/2}\leq C\left(\sum\limits_{m=0}^{+\infty}\left(\frac{\sqrt{40}b}{(1-\sqrt{\gamma})}\right)^m-1\right)=O(\delta_n)\,,
\end{eqnarray*}
provided that $1<b<\frac{\delta_n^{-1}(1-\sqrt{\gamma})}{\sqrt{40}}$ for some $n$ large enough. By Markov's inequality, it follows that $\theta_n\cip 0$, which in turn establishes \eqref{eq:prove1}.
}

\subsubsection{Proof of \eqref{eq:constant}}\label{sec:constant}
Set $$\eta_p=\frac{1}{2}\tr(\diag(\bR^{1/2}\tilde{\bS}\bR^{1/2})-\bI)^2=\frac{1}{2}\sum\limits_{i=1}^p\left(\frac{\tilde{\by}_i^\top\tilde{\by}_i}{n}-1\right)^2\,.$$ Denote the elements of the matrix $\bR^{1/2}=\{r_{(1/2),i,j}\}_{i,j=1}^p$ and let $\mathbf{r}_{(1/2),i}$ be the $i$th row of $\bR^{1/2}$, while the vector $\tilde{\by}_i$ is the $i$th row of the matrix $\bR^{1/2}\bX$. In order to proceed we need to rewrite the vector $\tilde{\by}_i$ in the following way
\begin{eqnarray}
  \tilde{\by}_i=\bX^\top\mathbf{r}_{(1/2),i}
\end{eqnarray}
and, thus,
\begin{eqnarray}\label{eta_revised} \frac{\tilde{\by}_i^\top\tilde{\by}_i}{n}&=&\frac{1}{n}\mathbf{r}^\top_{(1/2),i}\bX\bX^\top\mathbf{r}_{(1/2),i}=\frac{1}{n}\sum\limits_{j=1}^n\bx_j^\top\mathbf{r}_{(1/2),i}\mathbf{r}^\top_{(1/2),i}\bx_j=\frac{1}{n}\sum\limits_{j=1}^n\bx_j^\top\bR_i\bx_j\,.
\end{eqnarray}
Here we denoted for simplicity $\bR_i=\mathbf{r}_{(1/2),i}\mathbf{r}^\top_{(1/2),i}$.
Note that from the construction the following elementary identities hold
\begin{eqnarray}\label{indent_r12}
&& \bR=\bR^{1/2}\bR^{1/2}=\sum\limits_{i=1}^p\mathbf{r}_{(1/2),i}\mathbf{r}^\top_{(1/2),i}=\sum\limits_{i=1}^p\bR_i,\\
&& \tr(\bR_i)=\tr(\bR^2_i)=\mathbf{r}^\top_{(1/2),i}\mathbf{r}_{(1/2),i}=1\label{trace_Ri}\,.
% ~~\text{and}~~|\mathbf{r}^\top_{(1/2),i}\mathbf{r}_{(1/2),j}|<1\,.
\end{eqnarray}
Let us proceed to $\eta_p$, using \eqref{eta_revised}, \eqref{trace_Ri} and Lemma \ref{lem:quf} it holds
\begin{eqnarray*} \E(\eta_p)&=&\frac{1}{2}\sum\limits_{i=1}^p\E\left(\frac{\tilde{\by}^\top_i\tilde{\by}_i}{n}-1\right)^2= \frac{1}{2}\sum\limits_{i=1}^p\Var\left(\frac{\tilde{\by}^\top_i\tilde{\by}_i}{n}\right)= \frac{1}{2}\sum\limits_{i=1}^p\Var\left(\frac{1}{n}\sum\limits_{j=1}^n\bx_j^\top\bR_i\bx_j\right)\nonumber\\
                             &=& \frac{1}{2n^2}\sum\limits_{i=1}^p\sum\limits_{j=1}^n\Var(\bx_j^\top\bR_i\bx_j)+\frac{1}{n^2}\sum\limits_{i=1}^p\sum\limits_{k>l}\underbrace{\mathbbm{Cov}(\bx_k^\top\bR_i\bx_k, \bx_l^\top\bR_i\bx_l)}_{=0}\nonumber\\
                             &=& \frac{1}{2n^2}\sum\limits_{i=1}^p\sum\limits_{j=1}^n(\mathbbm{E}(x^4_{11})-3)\sum\limits_{k=1}^p(\bR_i)^2_{kk}+ \frac{1}{n^2}\sum\limits_{i=1}^p\sum\limits_{j=1}^n\tr(\bR^2_i)\\
                                             &=& \frac{1}{2}\frac{p}{n}(\mathbbm{E}(x^4_{11})-3)\frac{1}{p}\sum\limits_{i=1}^p\sum\limits_{k=1}^p(\mathbf{r}_{(1/2),i}\mathbf{r}^\top_{(1/2),i})^2_{kk}+\frac{p}{n}\\ &=&\frac{1}{2}\frac{p}{n}(\mathbbm{E}(x^4_{11})-3)\frac{1}{p}\sum\limits_{i=1}^p\sum\limits_{k=1}^pr^4_{(1/2),ik}+\frac{p}{n}\,.
\end{eqnarray*}
Denoting $\bi_p=(1,\underset{\text{$p$ times}}\ldots, 1)^\top$ and using the properties of the Hadamard product one can simplify the sum $\sum\limits_{i=1}^p\sum\limits_{k=1}^pr^4_{(1/2),ik}$ as follows
\begin{eqnarray*}
  \sum\limits_{i=1}^p\sum\limits_{k=1}^pr^4_{(1/2),i,k}&=&\bi_p^\top \Big( \bR^{1/2}\underset{\text{4 times}}{\circ\cdots\circ}\bR^{1/2}\Big)\bi_p=\tr\left(\bR^{1/2}\underset{\text{4 times}}{\circ\cdots\circ}\bR^{1/2}\bi_p\bi_p^\top \right)\\
                                                      &=& \tr\left([\bR^{1/2}\circ\bR^{1/2}][(\bR^{1/2}\circ\bR^{1/2})\circ\bi_p\bi_p^\top] \right)\\
  &=& \tr\left([\bR^{1/2}\circ\bR^{1/2}]^2\right)\,,
\end{eqnarray*}
which (recalling the definition of $C_{\bR^{1/2}}$) implies
\begin{eqnarray}\label{eq:expeta}
\mathbbm{E}(\eta_p)=\frac{1}{2}\frac{p}{n}(\mathbbm{E}(x^4_{11})-3)\frac{1}{p}\text{tr}\left([\bR^{1/2}\circ\bR^{1/2}]^2\right)+\frac{p}{n}=\frac{1}{2}\frac{p}{n}(\mathbbm{E}(x^4_{11})-3)C_{\bR^{1/2}}+\frac{p}{n}\,.
\end{eqnarray}

To prove \eqref{eq:constant}, we will show that the variance of $\eta_p$ converges to zero as $n \to \infty$. We have
\begin{equation}\label{eq:vareta}
4 \Var(\eta_p) = 4 \mathbbm{E}(\eta_p^2) - 4 (\mathbbm{E}(\eta_p))^2
\end{equation}
and $\mathbbm{E}(\eta_p)$ is given in \eqref{eq:expeta}.
For $1\le i\le p$ and $1\le t\le n$ set $Q_{it}= Q_{it}^{(n)}=\bx_t^\top\bR_i\bx_t-1$
and observe that $\E[Q_{it}]=0$. Since $\bx_1, \ldots, \bx_n$ are independent, it follows that $\E[Q_{it_1}Q_{jt_2}]=0$ whenever $t_1\neq t_2$.
By definition of $\eta_p$ and \eqref{eta_revised}, we have
\begin{equation*}
\begin{split}
2\eta_p &= \sum_{i=1}^p\left(\frac{\tilde{\by}_i^\top\tilde{\by}_i}{n}-1\right)^2
= \frac{1}{n^2} \sum_{i=1}^p \left( \sum_{t=1}^n (\bx_t^\top\bR_i\bx_t-1) \right)^2\\
&= \frac{1}{n^2} \sum_{i=1}^p \left( \sum_{t=1}^n Q_{it} \right)^2\,.
\end{split}
\end{equation*}
Therefore we get
\begin{equation}\label{eq:34}
\begin{split}
4 n^4\, \mathbbm{E}(\eta_p^2) &= \sum_{i,j=1}^p \sum_{t_1,\ldots, t_4=1}^n \E[Q_{it_1}Q_{it_2}Q_{jt_3}Q_{jt_4}]\\
&= \sum_{i,j=1}^p \Big( \sum_{t=1}^n \E[Q_{it}^2 Q_{jt}^2]
+ \sum_{t_1\neq t_2} \E[Q_{it_1}^2]\, \E[Q_{jt_2}^2]
+ 2 \sum_{t_1\neq t_2} \E[Q_{it_1} Q_{jt_1}] \,\E[Q_{it_2} Q_{jt_2}]\Big)\,.
\end{split}
\end{equation}

In view of $\E[Q_{it} Q_{jt}]= \E[\bx_t^\top\bR_i\bx_t \bx_t^\top\bR_j\bx_t] -1$ and $\tr(\bR_i)= \tr(\bR_i^2)=1$, an application of Lemma \ref{lem:quf} yields
\begin{equation}\label{eq:momq}
\E[Q_{it} Q_{jt}]= 2 \underbrace{\tr (\bR_i \bR_j)}_{\le \sqrt{\tr \bR_i^2 \tr \bR_j^2}} + (\E[x_{11}^4]-3) \underbrace{\tr (\bR_i \circ \bR_j)}_{\le \tr(\bR_i) \tr(\bR_j)} \le 2+ \E[x_{11}^4]\,.
\end{equation}
This implies that for a constant $C>0$ not depending on $n$ (which may change from one appearance to the next), we have
\begin{equation}\label{eq:segse}
\begin{split}
\sum_{i,j=1}^p \E[Q_{it_1} Q_{jt_1}] \,\E[Q_{it_2} Q_{jt_2}] &\le C \sum_{i,j=1}^p \E[Q_{it_1} Q_{jt_1}]\\
&\le C \sum_{i,j=1}^p \left( 2 \tr (\bR_i \bR_j) + (\E[x_{11}^4]-3) \tr (\bR_i \circ \bR_j) \right)\\
&\le C (\tr(\bR^2) + \tr(\bR \circ \bR)) \le C( p \norm{\bR} +p)\\
&\le C \,p\,.
\end{split}
\end{equation}
Here we used $\sum_i \bR_i = \bR$ and $\tr(\bR^2)\le \norm{\bR} \tr(\bR)$.

Next, set $\kappa= \E[x_{11}^4]-3$. By \eqref{eq:momq} and the definition of $C_{\bR^{1/2}}$, we have $\sum_{i=1}^p \E[Q_{it}^2]=p (2+ \kappa C_{\bR^{1/2}})$.
It easily follows that
\begin{equation*}
\begin{split}
\sum_{i,j=1}^p \E[Q_{it_1}^2]\, \E[Q_{jt_2}^2] &= \Big( \sum_{i=1}^p \E[Q_{it_1}^2]\Big)^2= p^2 (2+ \kappa C_{\bR^{1/2}})^2\,.
\end{split}
\end{equation*}
Hence, we deduce that
\begin{equation}\label{eq:setscdsw}
\begin{split}
\frac{1}{n^4} \sum_{t_1\neq t_2} \sum_{i,j=1}^p \E[Q_{it_1}^2]\, \E[Q_{jt_2}^2] - 4 (\mathbbm{E}(\tilde{\eta}_p))^2
&= \frac{n(n-1)p^2}{n^4} (2+ \kappa C_{\bR^{1/2}})^2 - \frac{p^2}{n^2} (2+ \kappa C_{\bR^{1/2}})^2\\
&= -\frac{p^2}{n^3} (2+ \kappa C_{\bR^{1/2}})^2
\end{split}
\end{equation}
A combination of \eqref{eq:vareta}, \eqref{eq:34}, \eqref{eq:segse} and \eqref{eq:setscdsw} shows that
\begin{equation}\label{eq:newvar}
4 \Var(\eta_p) = \frac{1}{n^4} \sum_{i,j=1}^p  \sum_{t=1}^n \E[Q_{it}^2 Q_{jt}^2] +O(1/n)\,,\qquad \nto\,.
\end{equation}
It remains to prove that
\begin{equation}\label{eq:left}
\frac{1}{n^4} \sum_{i,j=1}^p  \sum_{t=1}^n \E[Q_{it}^2 Q_{jt}^2]=o(1)\,.
\end{equation}
To this end we bound $\E[Q_{it}^2 Q_{jt}^2]$. Write $\bR_i=(R_{i,kl})$. Using the inequality $(a+b)^2\le 2 (a^2 +b^2)$ we get
\begin{equation*}
Q_{it}^2= \Big( \underbrace{\sum_{k=1}^p  R_{i,kk} (x_{kt}^2-1)}_{:=A_{it}} +\underbrace{\sum_{k\neq \ell}  R_{i,k\ell} x_{kt} x_{\ell t}}_{:=B_{it}}\Big)^2 \le 2 (A_{it}^2 +B_{it}^2)
\end{equation*}
and therefore
\begin{equation}\label{eq:boundq2}
\E[Q_{it}^2 Q_{jt}^2] \le 4 \left( \E[A_{it}^2 A_{jt}^2] +\E[B_{it}^2 B_{jt}^2] + \E[A_{it}^2 B_{jt}^2] +\E[A_{jt}^2 B_{it}^2]  \right)\,.
\end{equation}
We proceed by bounding the terms in \eqref{eq:boundq2}. Regarding the first term, a  direct calculation yields
\begin{equation}\label{eq:firstterm}
\E[A_{it}^2 A_{jt}^2]= \sum_{k=1}^p R_{i,kk}^2 R_{j,kk}^2 \E[(x_{kt}^2-1)^4] \le C \, \E[x_{11}^8] \sum_{k=1}^p R_{i,kk}^2 R_{j,kk}^2\le C \, n^2 \, \delta_n^4 \sum_{k=1}^p R_{i,kk}^2 R_{j,kk}^2\,.
\end{equation}
For the second term, an application of Lemma~\ref{lem:qmoment} and $\tr(\bR_i^2)=1$ yield
\begin{equation*}
(\E[B_{it}^4])^{1/4} = \Big( \E\Big[ \Big| \sum_{k\neq \ell}  R_{i,k\ell} x_{kt} x_{\ell t}\Big|^4\Big] \Big)^{1/4} \le C \, (\E[x_{11}^4])^{1/2} \Big(\sum_{k\neq \ell} R_{i,k\ell}^2\Big)^{1/2} \le C\,.
\end{equation*}
By Cauchy-Schwarz, we have
\begin{equation}\label{eq:secondterm}
\E[B_{it}^2 B_{jt}^2]\le \sqrt{\E[B_{it}^4] \, \E[B_{jt}^4]} \le C\,.
\end{equation}
Now we turn to the third and fourth terms. We have
\begin{equation*}
\E[A_{it}^2 B_{jt}^2] \le \sqrt{\E[A_{it}^4] \, \E[B_{jt}^4]}\le C n \, \delta_n^2 \Big(\sum_{k=1}^p R_{i,kk}^4\Big)^{1/2} \le C n \, \delta_n^2\,
\end{equation*}
from which it easily follows that
\begin{equation}\label{eq:thirdterm}
\E[A_{it}^2 B_{jt}^2] +\E[A_{jt}^2 B_{it}^2] \le C n \, \delta_n^2\,.
\end{equation}
Finally, we prove \eqref{eq:left}. Combining \eqref{eq:boundq2} and the inequalities for the individual terms (namely \eqref{eq:firstterm}, \eqref{eq:secondterm} and \eqref{eq:thirdterm}), we obtain
\begin{equation*}
\begin{split}
\frac{1}{n^4} \sum_{i,j=1}^p  \sum_{t=1}^n \E[Q_{it}^2 Q_{jt}^2]
&\le \frac{C}{n^3} \sum_{i,j=1}^p \left( n^2 \, \delta_n^4 \sum_{k=1}^p R_{i,kk}^2 R_{j,kk}^2 +1 +n \delta_n^2 \right) \\
&= \frac{C \, \delta_n^4}{n} \sum_{k=1}^p \left( \sum_{i=1}^p R_{i,kk}^2 \right)^2 + O(1/n)+O(\delta_n^2)\,.
\end{split}
\end{equation*}
Since $\sum_{i=1}^p \bR_i=\bR$ and $\diag(\bR)=\bI$, we have $\sum_{i=1}^p R_{i,kk}^2 \le \sum_{i=1}^p R_{i,kk} =1$. This completes the proof of \eqref{eq:left}.

\subsubsection{The case $p/n\to1$.}\label{sec:pn1}

{
One important difference in the case $p/n\to 1$ is that the  variance $\sigma_n^2$ tends to infinity. More precisely, it holds that
\begin{equation}\label{eq:sigma1}
\sigma_n^2=  -2\log\left(1-\frac{p-1}{n}\right)-2\frac{p}{n}+2\frac{p}{n}\tr(\bR-\bI)^2/p \sim -2\log\left(1-\frac{p-1}{n}\right)\to \infty\,.
\end{equation}
In Section \ref{sec:jointnormal}, we have shown that the last summand in \eqref{decomp} is $O_{\P}(1)$, which gives us $\log \det (\hat{\bR}) =\log \det (\bR^{1/2}\tilde{\bS}\bR^{1/2}) +O_{\P}(1)$. 
Using that $\sigma_n^2\to \infty$ and Slutsky's lemma, we deduce that the CLT in \eqref{CLT} follows from 
\begin{eqnarray}\label{CLT1}
    \frac{\log \det (\tilde{\bS})-\ov\mu_n}{\ov\sigma_n}\overset{d}{\longrightarrow}\mathcal{N}(0, 1),\qquad \nto\,,
  \end{eqnarray}
	where
\begin{align*}
		\ov\mu_n&= \left(p-n+\frac{1}{2}\right)\log\left(1-\frac{p-1}{n}\right)-p \quad \text{ and } \quad
    \ov\sigma_n^2=  -2\log\left(1-\frac{p-1}{n}\right).
  \end{align*}

It remains to prove \eqref{CLT1}. To this end, we refine the result of \cite{wang2018} using the unified expression. So, due to \cite{wang2018}  we get that as $p/n\to 1$ and $p<n$
  \begin{eqnarray}\label{cltwang1}
    \frac{\log \det(\tilde{\bS})-\sum\limits_{i=1}^p\log(1-i/n)}{ \ov\sigma_n}\overset{d}{\longrightarrow}\mathcal{N}(0, 1)\,
  \end{eqnarray}
  and otherwise if $p=n$ we have
  \begin{eqnarray}\label{cltwang2}
     \frac{\log \det(\tilde{\bS})+n\log n-\log(n-1)!}{\sqrt{2\log n}}\overset{d}{\longrightarrow}\mathcal{N}(0, 1)\,.
  \end{eqnarray}
	
We start with the case $p=n$, where $\ov \sigma_n^2=2 \log n$.  Using the Stirling formula $ \log(n!)=n\log(n)-n+(1/2)\log(2\pi n)+O(1/n)$, it is straightforward to show that
\begin{equation*}
\frac{n\log n-\log(n-1)!}{\sqrt{2\log n}}= \frac{\tfrac{1}{2}\log n +n}{\sqrt{2\log n}}+o(1)= \frac{\ov\mu_n}{\ov\sigma_n}+o(1).
\end{equation*}
In view of \eqref{cltwang2}, this establishes \eqref{CLT1} in the case $p=n$.
	
Next, we turn to the case $p/n\to1$ and $p<n$. Taking the logarithm on both sides of the identity
	\begin{eqnarray*}%\label{equal_prob}
      \prod\limits_{i=1}^{p}\left(1-\frac{i}{n}\right)=\frac{n!(1-p/n)}{(n-p)!\, n^p}=\frac{n!(n-p)}{n(n-p)!\, n^p}=\frac{(n-1)!}{(n-p-1)!\, n^p}\,,
  \end{eqnarray*}
we get 
  \begin{eqnarray}\label{simple}
    \sum\limits_{i=1}^p\log(1-i/n)=\log(n-1)!-p\log n-\log(n-p-1)! \,, \qquad p<n\,.
  \end{eqnarray}
Using \eqref{simple} one can rewrite the CLT in \eqref{cltwang2} for $\log \det (\tilde{\bS})$ in the case $p<n$ as 
  \begin{eqnarray*}%\label{unif_clt}
    \frac{\log \det(\tilde{\bS})-\log(n-1)!+p\log n+\log(n-p-1)!}{\ov\sigma_n}\overset{d}{\longrightarrow}\mathcal{N}(0, 1)\,.
  \end{eqnarray*}
  Now we apply the Stirling formula to approximate the centering terms for $p\le n-2$ and $n\to\infty$ as follows 
  \begin{eqnarray*}
   && p\log(n)+ \log(n-p-1)!- \log(n-1)!\\
    &=&p\log(n)+(n-p-1)\log(n-p-1)-(n-p-1)+1/2\log\left(2\pi(n-p-1)\right)\\
                            && -(n-1)\log(n-1)+(n-1)-1/2\log\left(2\pi(n-1)\right)+O(1)\\
                                         &=& (n-1)\log\left(\frac{n-p-1}{n-1}\right)-p\log\left(\frac{n-p-1}{n}\right)+\tfrac{1}{2}\log\left(\frac{n-p-1}{n-1}\right)+p+O(1)\\
   &=& -\big(p-n+\tfrac 12\big)\log\left(1-\frac{p+1}{n}\right)+p+O(1)\\
	&=& -\big(p-n+\tfrac 12\big)\log\left(1-\frac{p-1}{n}\right)+p+O(1)\,.
  \end{eqnarray*}
Since $\ov\sigma_n\to \infty$ we conclude that for $p\le n-2$
\begin{equation}\label{eq:sdde}
\frac{p\log(n)+ \log(n-p-1)!- \log(n-1)!}{\ov\sigma_n}= -\frac{\ov\mu_n}{\ov\sigma_n} +o(1)\,.
\end{equation}
A similar argument shows that \eqref{eq:sdde} also holds if $p=n-1$. In view of \eqref{cltwang1}, this establishes \eqref{CLT1} in the case $p/n\to1$ and $p<n$.
}

\subsection{Proof of Theorem \ref{sp}}\label{proofs:sp}

{
The proof is very similar to the proof of the main theorem, thus, we will skip some parts for the sake of simplicity and brevity.\\

The centered sample correlation matrix $ \hat \bR_c$ can be rewritten in terms of the non-centered sample covariance matrix $\bS$ in the following way
\begin{eqnarray}
  \hat \bR_c= \diag(\bS-\bar{\by}\bar{\by}^\top)^{-1/2} \,\left(\bS -  \bar{\by}\bar{\by}^\top\right)\diag(\bS-\bar{\by}\bar{\by}^\top)^{-1/2} \,,
\end{eqnarray}
where $\bar{\by}=1/n\bSigma^{1/2}\bX\bi=\bSigma^{1/2}\bar{\bx}$. Now using some matrix calculus we get
\begin{eqnarray*}
 \log \det (\hat \bR_c)&=& \log \det (\bS-\bar{\by}\bar{\by}^\top)-\log \det \left(\diag(\bS-\bar{\by}\bar{\by}^\top)\right)\\
&=& \log \det (\bR^{1/2}(\tilde{\bS}-\bar{\bx}\bar{\bx}^\top)\bR^{1/2}) - \log \det (\diag(\bR^{1/2}(\tilde{\bS}-\bar{\bx}\bar{\bx}^\top)\bR^{1/2}))\,.
\end{eqnarray*}

Again, denote the vector $\tilde{\by}_i$ as the $i$th row of the matrix $\bR^{1/2}\bX$, it holds
\begin{eqnarray*}
 \log \det (\diag(\bR^{1/2}(\tilde{\bS}-\bar{\bx}\bar{\bx}^\top)\bR^{1/2}))&=&\log \det \left(\frac{1}{n}\diag(\bR^{1/2}\bX(\bI-1/n\bi\bi^\top)\bX^\top\bR^{1/2})\right)\\
& =&\sum\limits_{i=1}^p\log(\tilde{\by}^\top_{i}(\bI-1/n\bi\bi^\top)\tilde{\by}_i/n)\,.
\end{eqnarray*}
One can check by Lemma \ref{bai_silverstein} that 
$$\max_{i=1,\ldots,p} \left|\tilde{\by}^\top_{i}(\bI-1/n\bi\bi^\top)\tilde{\by}_i/n-\frac{n-1}{n}\right|\to 0$$ in probability.
Therefore, using Taylor expansion of $\log(\tilde{\by}^\top_{i}(\bI-1/n\bi\bi^\top)\tilde{\by}_i/n)$ around $\frac{n-1}{n}$ we get
\begin{eqnarray*}
  \sum\limits_{i=1}^p\log(\tilde{\by}^\top_{i}(\bI-1/n\bi\bi^\top)\tilde{\by}_i/n)&=& p\log(1-1/n)+
\sum\limits_{i=1}^p\left(\frac{n}{n-1}\tilde{\by}^\top_{i}(\bI-1/n\bi\bi^\top)\tilde{\by}_i/n-1\right)
\\
&-&\frac12 \sum\limits_{i=1}^p\left(\frac{n}{n-1}\tilde{\by}^\top_{i}(\bI-1/n\bi\bi^\top)\tilde{\by}_i/n-1\right)^2\\
&+&\sum\limits_{i=1}^p\sum\limits_{q=3}^{+\infty}\frac{(-1)^{q+1}}{q}\left(\frac{n}{n-1}\tilde{\by}^\top_{i}(\bI-1/n\bi\bi^\top)\tilde{\by}_i/n-1\right)^q\\
                %                                              &=&\text{tr}(\tilde{\bS}^*-\bI)-\frac{1}{2}\text{tr}(\tilde{\bS}^*-\bI)^2+\sum\limits_{i=1}^pO(|\tilde{\bx}^\top_{i}\tilde{\bx}_i/n-1|^3)\,,\\
                                                              &=&p\log(1-1/n)+\tr\left(\frac{n}{n-1}\bR^{1/2}(\tilde{\bS}-\bar{\bx}\bar{\bx}^\top)\bR^{1/2}-\bI\right)\\
&-&\tilde{\eta}_p+\tilde{\xi}_{n,p}\,,
\end{eqnarray*}
where $\tilde{\eta}_p$ and $\tilde{\xi}_{n,p}$ are the expressions from the second and the third lines, respectively.
Altogether this gives
\begin{eqnarray*}
 \log \det (\hat{\bR}_c)&=&  \log \det (\bR^{1/2}(\tilde{\bS}-\bar{\bx}\bar{\bx}^\top)\bR^{1/2})-\tr\left(\frac{n}{n-1}\bR^{1/2}(\tilde{\bS}-\bar{\bx}\bar{\bx}^\top)\bR^{1/2}-\bI\right)\\
&& -p\log(1-1/n)+\frac{1}{2}\text{tr}\left(\text{diag}(\frac{n}{n-1}\bR^{1/2}(\tilde{\bS}-\bar{\bx}\bar{\bx}^\top)\bR^{1/2})-\bI\right)^2-\tilde{\xi}_{n,p}\,.
\end{eqnarray*}
Now we use the matrix determinant lemma to get
\begin{eqnarray*}
\log \det (\bR^{1/2}(\tilde{\bS}-\bar{\bx}\bar{\bx}^\top)\bR^{1/2}) = \log\left(1- \bar{\bx}^\top\tilde{\bS}^{-1}\bar{\bx}\right)+\log\det (\bR^{1/2}\tilde{\bS}\bR^{1/2})
\end{eqnarray*}
 and after some simplifications $\log\det(\hat{\bR}_c)$ can be rewritten as
  \begin{eqnarray}\label{asymp_log2}
    \log \det (\hat{\bR}_c)&=&  \log\left(1- \bar{\bx}^\top\tilde{\bS}^{-1}\bar{\bx}\right)+\bar{\bx}^\top\bR\bar{\bx}-p\log\left(1-\frac{1}{n}\right)-\frac{1}{n-1}\tr\left(\bR^{1/2}(\tilde{\bS}-\bar{\bx}\bar{\bx}^\top)\bR^{1/2}\right)\nonumber\\
&& + \log\det (\bR^{1/2}\tilde{\bS}\bR^{1/2})-\text{tr}(\bR^{1/2}\tilde{\bS}\bR^{1/2}-\bI) +\tilde{\eta}_p-\tilde{\xi}_{n,p}\,.
  \end{eqnarray}
First, we show that
$\tilde{\xi}_{n,p} \overset{\mathbbm{P}}{\longrightarrow} 0$.
Here we proceed similarly as in the proof of Theorem \ref{thm:main} by replacing the term $\tilde{\by}^\top_{i}\tilde{\by}_i/n-1$ by $\tilde{\by}^\top_{i}(\bI-1/n\bi\bi^\top)\tilde{\by}_i/n-(n-1)/n$. It can be seen that all of the arguments apply in the same way and the norm of the matrix $(\bI-1/n\bi\bi^\top)$ is obviously uniformly bounded. The next term we consider is
{
\begin{eqnarray}\label{tilde_eta}
 \tilde{\eta}_p&=&\frac12 \sum\limits_{i=1}^p\left(\frac{n}{n-1}\tilde{\by}^\top_{i}(\bI-1/n\bi\bi^\top)\tilde{\by}_i/n-1\right)^2= \frac{1}{2}\text{tr}(\text{diag}(\frac{n}{n-1}\bR^{1/2}(\tilde{\bS}-\bar{\bx}\bar{\bx}^\top)\bR^{1/2})-\bI)^2\nonumber\\
&=&\frac{1}{2}\text{tr}\left(\frac{1}{n-1}\text{diag}(\bR^{1/2}\bX(\bI-1/n\bi\bi^\top)\bX^\top\bR^{1/2})-\bI\right)^2\nonumber\\
&=& \frac{1}{2} \sum\limits_{i=1}^p\left(\frac{\tilde{\by}^\top_i(\bI-1/n\bi\bi^\top)\tilde{\by}_i}{n-1}-1  \right)^2=\frac{1}{2}\sum\limits_{i=1}^p\left(\frac{\tilde{\by}^\top_i\tilde{\by}_i}{n-1}-1 \right)^2-\sum\limits_{i=1}^p\left(\frac{\tilde{\by}^\top_i\tilde{\by}_i}{n-1}-1 \right)\bar{\bx}^\top\bR_i\bar{\bx}\frac{n}{n-1}\nonumber\\
&&+\frac{1}{2}\sum\limits_{i=1}^p\frac{n^2}{(n-1)^2}(\bar{\bx}^\top\bR_i\bar{\bx})^2\nonumber\\
&=&\eta_p-\sum\limits_{i=1}^p\left(\frac{\tilde{\by}^\top_i\tilde{\by}_i}{n-1}-1 \right)\bar{\bx}^\top\bR_i\bar{\bx}+\frac{1}{2}\sum\limits_{i=1}^p\frac{n^2}{(n-1)^2}(\bar{\bx}^\top\bR_i\bar{\bx})^2+o(1)\,,
\end{eqnarray}
where we have used $\tilde{\by}_i=\bX^\top\br_{(1/2),i}$ together with
\begin{eqnarray*}
  \frac{\tilde{\by}^\top_i(\bI-1/n\bi\bi^\top)\tilde{\by}_i}{n-1}=\frac{1}{n-1}\br_{(1/2),i}^\top\bX(\bI-1/n\bi\bi^\top)\bX^\top\br_{(1/2),i}=\frac{1}{n-1}\sum\limits_{j=1}^n\bx_j^\top\bR_i\bx_j-\frac{n}{n-1}\bar{\bx}^\top\bR_i\bar{\bx}\,
\end{eqnarray*}
with $\bR_i=\br_{(1/2),i}\br_{(1/2),i}^\top$. Moreover, it is easy to check that
\begin{eqnarray*}
\frac{1}{2} \sum\limits_{i=1}^p\left(\frac{\tilde{\by}^\top_i(\bI-1/n\bi\bi^\top)\tilde{\by}_i}{n-1}-1  \right)^2= \eta_p+o_{\mathbbm{P}}(1)\,.
\end{eqnarray*}
Now, from the proof of Theorem \ref{thm:main} we know that $\eta_p$ is asymptotically equivalent to $\frac{p}{2n}(\mathbbm{E}|x_{11}|^4-3)C_{\bR^{1/2}}+\frac{p}{n}$. Regarding the other terms in $\tilde{\eta}_p$ we will show that they converge to zero in probability. It has to be noted that this will follow from the fact that for $n\to\infty$
$$  \max\limits_{1\leq i\leq p} \bar{\bx}^\top\bR_i\bar{\bx}= \max\limits_{1\leq i\leq p} (\bar{\bx}^\top\br_{(1/2),i})^2\to 0~~\text{in probability}\,,
$$
where the latter is the consequence of the union bound, Markov inequality and Lemma \ref{bai_silverstein} for some $q$ large enough.
 Thus, we get for the second term in \eqref{tilde_eta}
\begin{eqnarray*}
 \left| \sum\limits_{i=1}^p\left(\frac{\tilde{\by}^\top_i\tilde{\by}_i}{n-1}-1 \right)\bar{\bx}^\top\bR_i\bar{\bx}\right| &\leq& \max\limits_{1\leq i\leq p} \bar{\bx}^\top\bR_i\bar{\bx}\cdot \left|\sum\limits_{i=1}^p\left(\frac{\tilde{\by}^\top_i\tilde{\by}_i}{n-1}-1 \right)\right|\\
&\leq& \max\limits_{1\leq i\leq p} \bar{\bx}^\top\bR_i\bar{\bx}\cdot \left|\sum\limits_{i=1}^p\left(\frac{\tilde{\by}^\top_i\tilde{\by}_i}{n}-1 \right)\right|(1+o(1)) \\
&=&\max\limits_{1\leq i\leq p} \bar{\bx}^\top\bR_i\bar{\bx}\cdot\left| \tr\left(\bR^{1/2}\tilde{\bS}\bR^{1/2}-\bI\right)\right|=o_{\mathbbm{P}}(1)O(1)=o_{\mathbbm{P}}(1)\,,
\end{eqnarray*}
because using Theorem \ref{thm:main} we know that $\tr\left(\bR^{1/2}\tilde{\bS}\bR^{1/2}-\bI\right)$ is a linear spectral statistics which is asymptotically normal with zero mean and, thus, it is bounded. Regarding the third term in \eqref{tilde_eta} we receive
\begin{eqnarray*}
  \frac{1}{2}\frac{n^2}{(n-1)^2}\sum\limits_{i=1}^p(\bar{\bx}^\top\bR_i\bar{\bx})^2\leq\left(1+\frac{1}{n-1} \right)^2 \frac{1}{2}\max\limits_{1\leq i\leq p} \bar{\bx}^\top\bR_i\bar{\bx}\cdot\sum\limits_{i=1}^p\bar{\bx}^\top\bR_i\bar{\bx}=o_{\mathbbm{P}}(1) \cdot \bar{\bx}^\top\bR\bar{\bx}=o_{\mathbbm{P}}(1)\,,
\end{eqnarray*}
where the last line follows from the fact that $\E(\bar{\bx}^\top\bR\bar{\bx})=p/n$ and $\Var(\bar{\bx}^\top\bR\bar{\bx})$ is vanishing. 

Therefore, we get that $\tilde{\eta}_p$ and $\eta_p$ are asymptotically equivalent. Moreover, the term $$\log\det (\bR^{1/2}\tilde{\bS}\bR^{1/2})-\text{tr}(\bR^{1/2}\tilde{\bS}\bR^{1/2}-\bI)$$
from \eqref{asymp_log2} is exactly the same as in Theorem \ref{thm:main}. Furthermore, it is easy to see using triangle inequality that for $n\to\infty$ in probability
\begin{eqnarray*}
   && \left|\bar{\bx}^\top\bR\bar{\bx} - \frac{1}{n-1}\tr\left(\bR^{1/2}(\tilde{\bS}-\bar{\bx}\bar{\bx}^\top)\bR^{1/2}\right) \right|\\
	&& \leq \left|\bar{\bx}^\top\bR\bar{\bx}-p/n\right|
+\left|\frac{1}{n-1}\tr\left(\bR^{1/2}(\tilde{\bS}-\bar{\bx}\bar{\bx}^\top)\bR^{1/2}\right)-p/n\right| \to 0\,,
\end{eqnarray*}
where the latter follows from \cite[Lemma 4]{rubiomestre2011} and \cite[Theorem A.44]{baisil2010}. 

Note also that for $n\to\infty$
\begin{eqnarray}
  (n-p+1/2)\left(\log\left(1-\frac{p}{n}\right)- \log\left(1-\frac{p}{n-1}\right)\right)\sim \frac{p}{n},\label{taylor_log1}\\
-p\log(1-1/n)\sim \frac{p}{n}\label{taylor_log2}\,,
\end{eqnarray}
where we have used the Taylor expansion of the logarithm.
Thus, taking into account \eqref{taylor_log1}, \eqref{taylor_log2} and \eqref{eq:seww} it is enough to show under different regimes, namely $p/n\to\gamma<1$ and $p/n\to\gamma=1$ that the $\log\left(1- \bar{\bx}^\top\tilde{\bS}^{-1}\bar{\bx}\right)$ is asymptotically equivalent to $\log(1-\frac{p}{n})$.

\subsubsection{The case $p/n\to\gamma <1$.}
 Next we show that for $p/n\to\gamma\in (0, 1)$ as $n\to\infty$ holds
\begin{eqnarray}\label{enough}
\text{log}(1-\bar{\bx}^\top\tilde{\bS}^{-1}\bar{\bx}) \overset{a.s.}{\longrightarrow} \log(1-\gamma)\,.
\end{eqnarray}
 Writing $\gamma=p/n$ and using \eqref{taylor_log1}, \eqref{taylor_log2} and \eqref{eq:seww} one can immediately see that this result corresponds to the substitution of $n$ by $n-1$ in the main formula of Theorem \ref{thm:main}. Let's consider now $\text{log}(1-\bar{\bx}^\top\tilde{\bS}^{-1}\bar{\bx})$. Using Woodbury matrix identity we get
\begin{eqnarray*}
\bar{\bx}^\top\tilde{\bS}^{-1}\bar{\bx}&=&\frac{1}{n^2}\bi^\top\bX^\top(1/n\bX\bX^\top)^{-1}\bX\bi=\lim\limits_{z\to0^+}\frac{1}{n^2}\bi^\top\bX^\top(1/n\bX\bX^\top-z\bI)^{-1}\bX\bi\\
&=& \lim\limits_{z\to0^+}\text{tr}\left[\frac{1}{\sqrt{n}}\bX^\top(1/n\bX\bX^\top-z\bI)^{-1}\frac{1}{\sqrt{n}}\bX\left(\frac{\bi\bi^\top}{n}\right)\right]\\
&=& \lim\limits_{z\to0^+}\text{tr}\left(\frac{\bi\bi^\top}{n}\right)+z\text{tr}\left[(1/n\bX^\top\bX-z\bI)^{-1}\left(\frac{\bi\bi^\top}{n}\right)\right]\\
&=&1+\lim\limits_{z\to0^+}z\text{tr}\left[(1/n\bX^\top\bX-z\bI)^{-1}\left(\frac{\bi\bi^\top}{n}\right)\right]\,.
\end{eqnarray*}
The application of \cite{baimpan2007}[see, Theorem 1] yields, as $\nto$,
\begin{eqnarray*}
z\text{tr}\left[(1/n\bX^\top\bX-z\bI)^{-1}\left(\frac{\bi\bi^\top}{n}\right)\right]\overset{a.s.}{\rightarrow} 
%\left(\frac{\bi\bi^\top}{n}\right)((\gamma-1)+\gamma zm(z))=
\gamma-1+\gamma zm(z)\,,
\end{eqnarray*}
where $m(z)$ is the limiting Stieltjes transform of $\tilde{\bS}$, namely the limit of $\frac{1}{p}\text{tr}\left( \tilde{\bS}-z\bI \right)^{-1}$. Now letting $z\to0^+$ and taking into account that $\lim\limits_{z\to0^+}m(z)<\infty$ we get
\begin{eqnarray}\label{qformx}
\bar{\bx}^\top\tilde{\bS}^{-1}\bar{\bx}\overset{a.s.}{\rightarrow} 1+\gamma-1=\gamma\,.
\end{eqnarray}
The result \eqref{enough} for $\gamma<1$ follows now from \eqref{qformx} and the continuous mapping theorem.

\subsubsection{The case $p/n\to1$.}

In this case we require $z$ to approach zero with a specific speed. We will need the following lemma, which follows from a combination of the proof of Theorem 1 by \cite{baimpan2007} and Lemma A.1 in \cite{bao2015}. This result is rather straightforward, thus, we omit here the details.
\begin{lemma}\label{rate_ST}
  Let $\bZ=(z_{ij})_{p\times n}$ be a random matrix with $n/2\leq p\leq n$ and $z_{ij}$ are i.i.d. with mean zero, variance one and $\mathbbm{E}|z_{11}|^4<\infty$. Denote $\alpha_n=n^{-1/6}$ and let $\ba_n\in \mathbbm{R}_1^n=\{\ba\in\mathbbm{R}^n: ||\ba||=1\}$. Then we have
  \begin{eqnarray*}
    \E\left(\ba_n^\top\left( \frac{1}{n}\bZ^\top\bZ+\alpha_n\bI \right)^{-1}\ba_n \right)&=&\frac{1-p/n}{\alpha_n}+ \frac{p}{n}m(\alpha_n)+O(n^{-1/6}),\\
\Var\left(\ba_n^\top\left( \frac{1}{n}\bZ^\top\bZ+\alpha_n\bI \right)^{-1}\ba_n \right)&=&O(n^{-1/3})\,
  \end{eqnarray*}
and
\begin{eqnarray*}
  m(\alpha_n)=2\left(\alpha_n+1-p/n+\sqrt{(\alpha_n+1-p/n)^2+\frac{4\alpha_n p}{n}} \right)^{-1}\,.
\end{eqnarray*}
\end{lemma}

First, denote $\varepsilon_n=1-p/n\to 0$ and let $\tilde{\bS}^{1/2}$ be a symmetric square root of $\tilde{\bS}$, then we have 
\begin{eqnarray*}
 \left| \bar{\bx}^\top\tilde{\bS}^{-1}\bar{\bx}-\bar{\bx}^\top(\tilde{\bS}+\alpha_n\bI)^{-1}\bar{\bx}\right|&=& \alpha_n\bar{\bx}^\top\tilde{\bS}^{-1}(\tilde{\bS}+\alpha_n\bI)^{-1}\bar{\bx}= \alpha_n \bar{\bx}^\top\tilde{\bS}^{-1/2}(\tilde{\bS}+\alpha_n\bI)^{-1}\tilde{\bS}^{-1/2}\bar{\bx}\\
&\leq& \frac{\alpha_n}{\lambda_{min}(\tilde{\bS}+\alpha_n\bI)}\bar{\bx}^\top\tilde{\bS}^{-1}\bar{\bx}\leq \frac{\alpha_n}{\lambda_{min}(\tilde{\bS})}\bar{\bx}^\top\tilde{\bS}^{-1}\bar{\bx}  \,,
\end{eqnarray*}
which reveals, e.g., that for $\varepsilon_n=O(n^{-1/6+\delta})$ with $\delta\in(0 , 1/6)$ it holds for sufficiently large $n$ that
\begin{eqnarray}\label{replaceS}
   \left| 1-\frac{\bar{\bx}^\top(\tilde{\bS}+\alpha_n\bI)^{-1}\bar{\bx}}{\bar{\bx}^\top\tilde{\bS}^{-1}\bar{\bx}}\right| \le \frac{\alpha_n}{(1-\sqrt{1-\varepsilon_n})^2}=o(1)\,
\end{eqnarray}
because the smallest eigenvalue of $\tilde{\bS}$ is asymptotically equivalent to $(1-\sqrt{p/n})^2$ as $\nto$.
Thus, we may replace $\bar{\bx}^\top\tilde{\bS}^{-1}\bar{\bx}$ by $\bar{\bx}^\top(\tilde{\bS}+\alpha_n\bI)^{-1}\bar{\bx}$ without altering the limit. Reusing the calculations for the case $p/n\to\gamma<1$ we get
\begin{eqnarray*}
  \bar{\bx}^\top(\tilde{\bS}+\alpha_n\bI)^{-1}\bar{\bx}=1-\alpha_n\text{tr}\left[(1/n\bX^\top\bX+\alpha_n\bI)^{-1}\left(\frac{\bi\bi^\top}{n}\right)\right]\,.
\end{eqnarray*}
Thus, we have
\begin{eqnarray}\label{diff_log}
  \log\left(1-\bar{\bx}^\top(\tilde{\bS}+\alpha_n\bI)^{-1}\bar{\bx}\right)-\log(\varepsilon_n)=\log\left(\frac{\alpha_n}{\varepsilon_n}\text{tr}\left[(1/n\bX^\top\bX+\alpha_n\bI)^{-1}\left(\frac{\bi\bi^\top}{n}\right)\right]\right)\,.
\end{eqnarray}
Using Lemma \ref{rate_ST} we get taking $\varepsilon_n=O(n^{-1/12})$
\begin{eqnarray*}
 && \mathbbm{E}\left(\frac{\alpha_n}{\varepsilon_n}\text{tr}\left[(1/n\bX^\top\bX+\alpha_n\bI)^{-1}\left(\frac{\bi\bi^\top}{n}\right)\right]\right)=1+p/n\frac{\alpha_n}{\varepsilon_n}m(-\alpha_n)+o(1)\\
&=& 1+2p/n\frac{\alpha_n}{\varepsilon_n}\left(\alpha_n+\varepsilon_n+\sqrt{(\alpha_n+\varepsilon_n)^2+\frac{4 p}{n}\alpha_n} \right)^{-1}+o(1)=1+O(1)\,.
\end{eqnarray*}
and
\begin{eqnarray*}
  \Var\left(\frac{\alpha_n}{\varepsilon_n}\text{tr}\left[(1/n\bX^\top\bX+\alpha_n\bI)^{-1}\left(\frac{\bi\bi^\top}{n}\right)\right]\right)=\frac{\alpha^2_n}{\varepsilon^2_n}O(n^{-1/3})=o(1)\,.
\end{eqnarray*}
Now using the Markov inequality we may conclude that \eqref{diff_log} is bounded in probability as $p/n\to1$ and because \eqref{diff_log} is additionally divided by $\sqrt{-2\log\varepsilon_n}\to\infty$ we get by \eqref{replaceS} and the triangle inequality that
\begin{eqnarray}
 \frac{1}{\sqrt{-2\log\varepsilon_n}}\log\left(\varepsilon^{-1}_n(1-\bar{\bx}^\top\tilde{\bS}^{-1}\bar{\bx})\right)=o_{\mathbbm{P}}(1)\,,
\end{eqnarray}
which finishes the proof of Theorem \ref{sp} in case $p/n\to1$.
}
}

\section{Auxiliary results} \label{sec:5}
\begin{lemma}[\cite{baisil2010}, Lemma 9.1]\label{bai_silverstein}
  Let $\bA=(a_{ij})$ be an $n\times n$ nonrandom matrix and $\bx=(x_1,\ldots, x_n)$ be a random vector of independent entries. Assume that $\mathbbm{E}(x_i)=0$, $\mathbbm{E}(x^2_i)=1$, $\sup_i\mathbbm{E}(x^4_i)= \nu_4<\infty$, and $|x_i|\leq\delta_n\sqrt{n}$. Then for any given $2\leq q\leq b\log(n\nu^{-1}_4\delta_n^4)$ and $b>1$, we have
  \begin{eqnarray*}
    \mathbbm{E}\left|\bx^\top\bA\bx-\tr(\bA) \right|^q\leq  \nu_4n^q(40b^2)^q ||\bA||^q \delta_n^{2q-4}n^{-1} \,.
  \end{eqnarray*}
\end{lemma}
{
We need the following lemma; see for example parts $b)$ and $d)$ of Theorem in \cite{wiens:1992}.
\begin{lemma}[Moments of quadratic forms]\label{lem:quf}
  Let $\bz=(Z_1,\ldots,Z_n)^\top$ be a random vector with i.i.d. entries, with $\E[Z_1]=0, \E[Z_1^2]=\nu_2, \E[Z_1^4]=\nu_4< \infty$, and let $\bA,\bB, \bC$ be real and symmetric $n\times n$ nonrandom matrices. Then
  \begin{equation*}
	\E[\bz^{\top} \bA \bz \cdot \bz^{\top} \bB \bz]=  \tr(\bA) \tr(\bB) +2 \tr(\bA \bB) + (\nu_4-3) \tr(\bA \circ \bB)\,,
  \end{equation*}
	where $\circ$ denotes the Hadamard product. As a special case we get the variance
        \begin{eqnarray*}
          \Var(\bz^\top\bA\bz)=2\tr(\bA^2)+(\nu_4-3)\tr(\bA \circ \bA)\,.
        \end{eqnarray*}
 If aditionally $\E[Z_1^6]=\nu_6<\infty$, one has
        \begin{eqnarray*}
          \E[\bz^{\top} \bA \bz\cdot \bz^{\top} \bB \bz \cdot \bz^{\top} \bC \bz]&=&\tr \bA \tr \bB \tr\bC+2\left( \tr\bA\cdot\tr(\bB\bC) + \tr\bB\cdot\tr(\bA\bC)+\tr\bC\cdot\tr(\bA\bB)\right) \\
&+&(\nu_4-3) \left( \tr\bA\cdot\tr(\bB\circ\bC) + \tr\bB\cdot\tr(\bA\circ\bC)+\tr\bC\cdot\tr(\bA\circ\bB) \right)\\
&+& 4(\nu_4-3)\left(\tr(\bA\cdot(\bB\circ\bC))+\tr(\bB\cdot(\bA\circ\bC))+\tr(\bC\cdot(\bA\circ\bB)) \right)\\
&+& (\nu_6-15\nu_4+30)\tr(\bA\circ\bB\circ\bC)+ 8\tr(\bA\bB\bC)\,.
        \end{eqnarray*}
Particularly,
\begin{equation*}
	\E[(\bz^{\top} \bA \bz -\E[\bz^{\top} \bA \bz])^3]=  8 \tr(\bA^3) +12 (\nu_4-3) \tr(\bA \circ \bA^2) + (\nu_6-15 \nu_4 +30) \tr (\bA \circ \bA \circ \bA)\,.
  \end{equation*}	
\end{lemma}
}

The next result is Lemma 7.10 in \cite{erdos:yau:2017}.
\begin{lemma}\label{lem:qmoment}
Let $X_1,\ldots,X_N$ be independent centered random variables and assume that
\begin{equation*}
(\E[|X_i|^q])^{1/q} \le \mu_q\,,\quad 1\le i\le N; q=2,3,\ldots
\end{equation*}
for some fixed constants $\mu_q$. Then we have for any deterministic complex numbers $a_{ij}, 1\le i,j\le N$ that
\begin{equation*}
\Big(\E\Big[ \Big|\sum_{i\neq j=1}^N a_{ij}X_i X_j \Big|^q \Big]\Big)^{1/q}\le C \, q\, \mu_q^2 \Big(\sum_{i\neq j=1}^N |a_{ij}|^2\Big)^{1/2}\,,\quad q=2,3,\ldots,
\end{equation*}
where the constant $C$ does not depend on $q$.
\end{lemma}

\subsection*{Acknowledgments}

J.Heiny's research was supported by the Deutsche Forschungsgemeinschaft (DFG) via RTG 2131 High-dimensional Phenomena in Probability -- Fluctuations and Discontinuity. The authors thank Gabriella F.~Nane for fruitful discussions.

\setlength{\bibsep}{4pt}
%\begin{small}
\bibliography{literature}
%\end{small}

\end{document}